\documentclass{agtart_a}
\pdfoutput=1
\usepackage{graphicx}

\title[On the existence of branched coverings between
surfaces]{On the existence of branched coverings between
surfaces\\ with prescribed branch data, I}

\author{Ekaterina Pervova}
\givenname{Ekaterina}
\surname{Pervova}
\address{Chelyabinsk State University\\
ul. Br. Kashirinykh, 129\\\newline
454021 Chelyabinsk, Russia}
\email{pervova@csu.ru}
\urladdr{}

\author{Carlo Petronio}
\givenname{Carlo}
\surname{Petronio}
\address{Dipartimento di Matematica Applicata\\
Universit\`a di Pisa\\\newline
Largo Bruno Pontecorvo, 1\\
56127 Pisa, Italy}
\email{petronio@dm.unipi.it}
\urladdr{}

\volumenumber{6}
\issuenumber{}
\publicationyear{2006}
\papernumber{69}
\startpage{1957}
\endpage{1985}

\doi{}
\MR{}
\Zbl{}


\received{20 January 2006}
\revised{14 September 2006}
\accepted{25 September 2006}
\published{14 November 2006}
\publishedonline{14 November 2006}
\proposed{}
\seconded{}
\corresponding{}
\editor{CPR}
\version{}

\arxivreference{math.GT/0508434}

\AtBeginDocument{\let\bar\wbar\let\tilde\wtilde\let\hat\what
\renewcommand{\hbar}{{\bar{h}}}}

\makeatletter
\def\cnewtheorem#1[#2]#3{\newtheorem{#1}{#3}[section]
\expandafter\let\csname c@#1\endcsname\c@lemma}
\makeatother

\theoremstyle{plain}
\newtheorem{lemma}{Lemma}[section]
\cnewtheorem{thm}[lemma]{Theorem}
\cnewtheorem{prop}[lemma]{Proposition}
\cnewtheorem{cor}[lemma]{Corollary}
\cnewtheorem{conj}[lemma]{Conjecture}

\theoremstyle{remark}
\cnewtheorem{rem}[lemma]{Remark}
\cnewtheorem{example}[lemma]{Example}
\cnewtheorem{defn}[lemma]{Definition}
\cnewtheorem{fact}[lemma]{Fact}
\cnewtheorem{ques}[lemma]{Question}

\newcommand{\dimo}[1]{\proof[Proof of~\fullref{#1}]}
\newcommand{\finedimo}{\endproof}

\newcommand\matT{{\mathbb{T}}}
\newcommand\matS{{\mathbb{S}}}
\newcommand\matP{{\mathbb{P}}}
\newcommand\matC{{\mathbb{C}}}

\newcommand\Sigmatil{{\widetilde\Sigma}}
\newcommand\ntil{{\wtilde n}}

\newfont{\Got}{eufm10 scaled 1200}
\newcommand{\permu}{{\hbox{\Got S}}}

\begin{document}

\begin{asciiabstract}
For the existence of a branched covering Sigma~ --> Sigma between
closed surfaces there are easy necessary conditions in terms of
chi(Sigma~), chi(Sigma), orientability, the total degree, and the
local degrees at the branching points. A classical problem dating back
to Hurwitz asks whether these conditions are also sufficient.  Thanks
to the work of many authors, the problem remains open only when Sigma
is the sphere, in which case exceptions to existence are known to
occur. In this paper we describe new infinite series of exceptions, in
particular previously unknown exceptions with Sigma~ not the sphere
and with more than three branching points. All our series come with
systematic explanations, based on several different techniques
(including dessins d'enfants and decomposability) that we exploit to
attack the problem, besides Hurwitz's classical technique based on
permutations.  Using decomposability we also establish an easy
existence result.
\end{asciiabstract}

\begin{htmlabstract}
<p class="noindent">
For the existence of a branched covering
&Sigma;<sup>~</sup>&rarr;&Sigma; between closed surfaces there are easy necessary
conditions in terms of &chi;(&Sigma;<sup>~</sup>), &chi;(&Sigma;), orientability,
the total degree, and the local degrees at the branching points. A
classical problem dating back to Hurwitz asks whether these conditions
are also sufficient.  Thanks to the work of many authors, the problem
remains open only when &Sigma; is the sphere, in which case exceptions
to existence are known to occur. In this paper we describe new
infinite series of exceptions, in particular previously unknown
exceptions with &Sigma;<sup>~</sup> not the sphere and with more than three
branching points. All our series come with systematic explanations,
based on several different techniques (including dessins d'enfants and
decomposability) that we exploit to attack the problem, besides
Hurwitz's classical technique based on permutations.  Using
decomposability we also establish an easy existence result.  </p>
\end{htmlabstract}

\begin{abstract}
\noindent For the existence of a branched covering
$\widetilde{\Sigma}\to\Sigma$ between closed surfaces there are easy
necessary conditions in terms of $\chi(\widetilde{\Sigma})$,
$\chi(\Sigma)$, orientability, the total degree, and the local degrees
at the branching points. A classical problem dating back to Hurwitz
asks whether these conditions are also sufficient.  Thanks to the work
of many authors, the problem remains open only when $\Sigma$ is the
sphere, in which case exceptions to existence are known to occur. In
this paper we describe new infinite series of exceptions, in
particular previously unknown exceptions with $\widetilde{\Sigma}$ not
the sphere and with more than three branching points. All our series
come with systematic explanations, based on several different
techniques (including dessins d'enfants and decomposability) that we
exploit to attack the problem, besides Hurwitz's classical technique
based on permutations.  Using decomposability we also establish an
easy existence result.
\end{abstract}

\maketitle

\section{Problem and new partial solutions}\label{new:results:section}
In this section we state the Hurwitz existence problem, we
outline its relevance to other areas of topology
and our motivations for picking it up,
and we state our new contributions towards its solution, also
explaining the techniques we have used to obtain them.
We address the reader to \fullref{review:section}
for an overview of the known results and techniques,
which will help putting our results into context.

\paragraph{Basic definitions}
A \emph{branched covering} is a map
$f\co \Sigmatil\to\Sigma$, where $\Sigmatil$ and $\Sigma$ are closed
connected surfaces and $f$ is locally modelled on maps of the form
$\matC\ni z\mapsto z^k\in\matC$ for some $k\geqslant1$. The
integer $k$ is called the \emph{local degree} at the point of
$\Sigmatil$ corresponding to $0$ in the source $\matC$. If $k>1$
then the point of $\Sigma$ corresponding to $0$ in the target
$\matC$ is called a \emph{branching point}. The branching points
are isolated, hence there are finitely many, say $n$, of them.
Removing the branching points in $\Sigma$ and all their pre-images
in $\Sigmatil$, the restriction of $f$ gives a genuine covering,
whose degree we will denote by $d$. If the $i$-th branching point
on $\Sigma$ has $m_i$ pre-images, the local degrees
$(d_{ij})_{j=1,\ldots,m_i}$ at these points give a partition of
$d$, namely $d_{ij}\geqslant 1$ and $\sum_{j=1}^{m_i}d_{ij}=d$. In
the sequel we will always assume that in a partition
$(d_1,\ldots,d_m)$ of $d$ we have $d_1\geqslant\ldots\geqslant
d_m$, which allows us to regard the partition as an array of integers
rather than an unordered set with repetitions.

\paragraph{Branch data}
Suppose we are given closed connected surfaces $\Sigmatil$ and
$\Sigma$, integers $n\geqslant 0$ and $d\geqslant 2$, and for
$i=1,\ldots,n$ a partition $(d_{ij})_{j=1,\ldots,m_i}$ of $d$. The
5--tuple $\big(\Sigmatil,\Sigma,n,d,(d_{ij})\big)$ will be called
the \emph{branch datum} of a candidate branched covering. To such
a datum we will always associate the integer $\ntil$ defined as
$m_1+\ldots+m_n$.

\paragraph{Compatibility}
We define a branch datum to be \emph{compatible} if the following
conditions hold:
\begin{enumerate}
\item\label{RHcond}
$\chi(\Sigmatil)-\ntil=d\cdot(\chi(\Sigma)-n)$;
\item\label{Pcond}
$n\cdot d-\ntil$ is even;
\item\label{OOcond} If $\Sigma$ is
orientable then $\Sigmatil$ is also orientable;
\item\label{NNcond} If $\Sigma$ is non-orientable and $d$ is odd
then $\Sigmatil$ is also non-orientable;
\item\label{ONcond} If
$\Sigma$ is non-orientable but $\Sigmatil$ is orientable then each
partition $(d_{ij})_{j=1,\ldots,m_i}$ of $d$ refines the partition
$(d/2,d/2)$.
\end{enumerate}
The meaning of Condition~\ref{ONcond} is that
$(d_{ij})_{j=1,\ldots,m_i}$ is obtained by juxtaposing two
partitions of $d/2$ and reordering. Note that $d$ is even by
Condition~\ref{NNcond}.

\paragraph{The problem}
It is not too difficult to show that if a branched covering
$\Sigmatil\to\Sigma$ exists then the corresponding branch datum,
with $n,d,(d_{ij}),\ntil$ defined as above, is compatible.
Conditions~\ref{RHcond},~\ref{OOcond}, and~\ref{NNcond} are
obvious, Condition~\ref{ONcond} follows from the fact that the
covering factors through the orientation covering of $\Sigma$, and
a short proof of Condition~\ref{Pcond} will be given for the sake
of completeness in \fullref{views:section}.

We will call \emph{Hurwitz existence problem} the
question of which compatible branch data are actually realized by
some branched covering. In the sequel we will always consider the
branch datum $\big(\Sigmatil,\Sigma,n,d,(d_{ij})\big)$, and the
corresponding $\ntil$, to be fixed. We will also mostly assume
that each partition $(d_{ij})_{j=1,\ldots,m_i}$ is different from
$(1,\ldots,1)$, for in this case we could just reduce $n$.

\begin{rem}\label{no:Sigmatil:rem}
A perhaps more traditional viewpoint is to consider only
$\Sigma$, $n$, $d$, $(d_{ij})$ to be given, and then determine the
corresponding $\Sigmatil$ using
Conditions~\ref{RHcond},~\ref{OOcond}, and~\ref{NNcond}. In this
context, these conditions are replaced by the requirement that
$\ntil+d\cdot(\chi(\Sigma)-n)$ should be at most $2$. However, for
non-orientable $\Sigma$ and even $d$, two possibilities exist for
$\Sigmatil$, so we prefer to stick to our datum which also
includes $\Sigmatil$.\end{rem}

Even if we have stated the problem in full generality, we will
confine ourselves in the rest of this section to the case
$\Sigma=\matS$, the sphere, because a full solution has been
obtained (in the affirmative) whenever $\chi(\Sigma)\leqslant 0$,
and the case $\Sigma=\matP$, the projective plane, reduces to the
case $\Sigma=\matS$. Among the many sources for this solution,
that we will discuss in \fullref{review:section}, we single
out here the fundamental contribution of Edmonds, Kulkarni and
Stong~\cite{EKS}, that we will frequently refer to, and the more
recent one by Bar\'anski~\cite{Baranski}. We also mention the very
interesting paper by Zheng~\cite{Zheng}, which introduces a new approach to the
problem and describes many experimental results.

\paragraph{Motivation}
Surfaces are central objects in mathematics.  They are interesting
on their own (being the subject matter of, for instance,
Teichm\"uller theory) and they are relevant to diverse fields such
as algebraic geometry, complex analysis, and three-dimensional
topology. Branched coverings between surfaces naturally occur
within all these fields of investigation, so the basic Hurwitz
existence problem stated
above can be viewed in our opinion as one of
the crucial ones in modern mathematics.

As discussed in \fullref{review:section},
the problem has indeed attracted enormous attention over about a century,
including that of outstanding mathematicians.
The general solution of the problem for $\Sigma=\matS$ is however still
missing, which suggests that the problem is actually rather hard.

Besides being intrinsically interesting and very difficult, the
existence problem for branched coverings with prescribed branch
data naturally emerges in several contexts. For instance, it is
relevant to the study of generating sets of surface groups. The
reason is that a branched covering between surfaces naturally
induces an orbifold-covering between 2--orbifolds (see
\fullref{views:section} below), and coverings of the latter
type correspond to subgroups of orbifold fundamental groups.
Therefore existence of a branched covering matching a given datum
is equivalent to existence in a certain Fuchsian group of a
subgroup with given signature (Singerman~\cite{Singerman}), which can be
applied to deciding whether an arbitrary generating set of a
surface group is Nielsen-equivalent to the standard generating
set.

There are enumerative aspects of the Hurwitz problem which are not
directly faced in the present paper but have obvious connections with
our topic, and these aspects are relevant to the Gromov--Witten theory
of algebraic curves. The Hurwitz number associated to a given datum is
the number of equivalences classes of coverings realizing this datum
and, as indicated in Okounkov and Pandharipande~\cite{OkPand}, the
stationary Gromov--Witten invariant of a curve is equal to the sum of
the Hurwitz numbers associated to certain branch data, determined
through a specific correspondence between branch data and descendants
in Gromov--Witten theory.

Our own motivation for picking up the existence problem was to
investigate the behaviour of Matveev's complexity~\cite{Matveev} for
3--manifolds under finite covering. Conjectural formulae for
complexity have been given by Martelli and Petronio~\cite{MaPeGD} for
Seifert 3--manifolds, ie, for circle fibrations over
2--orbifolds, an important class of 3--manifolds which has been
classified for a long time (see eg Matveev and
Fomenko~\cite{MaFo}). It turns out that a finite covering between
Seifert 3--manifolds induces a covering between the corresponding base
2--orbifolds, and such a covering can be interpreted as a branched
covering between surfaces, as already mentioned above.  Therefore
understanding branched coverings between surfaces can be viewed as a
first necessary step towards the analysis of coverings between Seifert
3--manifolds and hence of Matveev's complexity under such coverings.

\paragraph{New results}
A branch datum will be called \emph{exceptional} if it is compatible
but not realizable by any branched covering. As already mentioned,
for $\Sigma=\matS$ exceptional data are known to exist. All examples
discussed in the literature refer to the case where $d$ is
non-prime, $n$ is $3$, and $\Sigmatil$ is also
equal to $\matS$.

The main results obtained in the paper are listed below. They improve
the understanding of exceptional data in that they place most of the
exceptions occurring for non-prime $d\leqslant 22$ and
$\Sigmatil=\Sigma=\matS$ (as described by Zheng
in~\cite{Zheng,Zheng-bis}) within infinite series of exceptions, thus
providing some sort of explanation for their emergence, and they show
that such systematic exceptions occur also for $\Sigmatil\neq\matS$.
Moreover, Theorems~\ref{excep:by:fixpoints:thm}
and~\ref{even-deg_exceptions:thm} below show that there are several
exceptional series of data with $n>3$ branching points.

A result of rather different nature is given by
\fullref{divisible:exist:thm}, which provides a simple condition for
realizability. More such conditions, based on a generalization of the
results of Bar\'anski~\cite{Baranski}, will be described
in~\cite{partII}.

\begin{prop}\label{non-realiz:53:prop}
Let $d\geqslant 8$ be even and consider compatible branch data of
the form $\big(\Sigmatil,\matS,3,d,(2,\ldots,2),
(5,3,2,\ldots,2), (d_{3j})_{j=1,\ldots,m_3}\big)$.
\begin{itemize}
\item If $\Sigmatil=\matT$, the torus, whence $m_3=2$, the datum is
realizable if and only if $(d_{31},d_{32})\neq(d/2,d/2)$; \item If
$\Sigmatil=\matS$, whence $m_3=4$, the datum is realizable if and
only if $(d_{3j})_{j=1,\ldots,4}$ does not have the form
$(k,k,d/2-k,d/2-k)$ for some $k>0$, or\break $(d/2,d/6,d/6,d/6)$ for $d$
a multiple of $6$.
\end{itemize}
\end{prop}

\begin{prop}\label{23:nonex:prop}
Let $\big(\matS,\matS,3,d,(d_{ij})\big)$ be a compatible branch
datum with even $d$ and $(d_{1j})=(2,\ldots,2)$. If
$(d_{2j})=(3,3,2,\ldots,2)$ or $(d_{2j})=(3,2,\ldots,2,1)$ then
the datum is realizable if and only if $d_{31}\neq d/2$.
\end{prop}

\begin{thm}\label{excep:by:fixpoints:thm}
Suppose that $d$ and all $d_{ij}$ for $i=1,2$ are multiples of
some $k$ with $1<k<d$. If the branch datum
$\big(\matS,\matS,n,d,(d_{ij})\big)$ is realizable then
$d_{ij}\leqslant d/k$ for $i=3,\ldots,n$ and for all $j$.
\end{thm}

\begin{thm}\label{even-deg_exceptions:thm}
Suppose that $d$ and all $d_{ij}$ for $i=1,2$ are even. If the
branch datum $\big(\matS,\matS,n,d,(d_{ij})\big)$ is realizable
then $(d_{ij})$ refines the partition $(d/2,d/2)$ for
$i=3,\ldots,n$.
\end{thm}

\begin{cor}\label{mixed:deg:excep:cor}
Suppose that $d$ is a multiple of $2k$ for some $k$ with
$1<k<d/2$, that all $d_{1j}$ are multiples of $k$, and that all
$d_{ij}$ for $i=2,3$ are even. If the branch datum
$\big(\matS,\matS,n,d,(d_{ij})\big)$ is realizable then
$d_{ij}\leqslant d/k$ for $i=2,3$ and $d_{ij}\leqslant d/2k$ for
$i=4,\ldots,n$ and for all $j$.
\end{cor}

The last three criteria are especially expected to cover very many
exceptional branch data. For example they imply that all the
following series of data are exceptional:
\begin{multline*}
\big(\matS,\matS,d/k+1,d,(k,\ldots,k),(k,\ldots,k),(d/k+1,1,\ldots,1),\\
\shoveright{(2,1,\ldots,1),\ldots,(2,1,\ldots,1) \big),\quad k|d,}\\
\shoveleft{\big(\matS,\matS,d/2+1,d,(2,\ldots,2),(2,\ldots,2),(2,\ldots,2),} \\
\shoveright{(2,1,\ldots,1),\ldots,(2,1,\ldots,1) \big),\quad
d\equiv 2\mod 4,}\\
\shoveleft{\big(\matS,\matS,d/2k+2,d,(k,\ldots,k),(2,\ldots,2),(2,\ldots,2),(d/2k+1,1,\ldots,1),}\\
(2,1,\ldots,1), \ldots,(2,1,\ldots,1) \big),\quad
2k|d.\end{multline*}
We note that the exceptionality of the first
of these series was already conjectured in general and proved for
$d\leqslant 20$ by Zheng in~\cite[Conjecture 16]{Zheng}.

\begin{thm}\label{divisible:exist:thm} Let
$\big(\Sigmatil,\matS,3,d,(d_{ij})\big)$ be a compatible branch
datum. Let $p\geqslant 3$ be odd and suppose that all $d_{ij}$ are
divisible by $p$. Then the datum is realizable.
\end{thm}

\paragraph{Techniques}
Various equivalent ways of formulating the Hurwitz existence
problem, and techniques to attack it, were developed over the
time. We only mention here that the main classical tool, which
goes back to Hurwitz himself, is a reformulation of the problem in
terms of permutations. The main techniques we employ are as
follows:

\begin{itemize}
\item \textbf{Dessins d'enfants}\qua This is a classical notion due to
Grothendieck, introduced within the study of algebraic maps between
Riemann surfaces.  This topic is tightly related to our existence
problem, but dessins d'enfants were never employed directly before to
attack the problem itself, and our adjustment of the notion proved
rather fruitful, leading to Propositions~\ref{non-realiz:53:prop}
and~\ref{23:nonex:prop}, and to
\fullref{excep:by:fixpoints:thm}.

\item
\textbf{Decomposability}\qua The second tool we use is based on the idea of
expressing a covering as a composition of two non-trivial ones, and
the corresponding idea of finding certain ``block decompositions''
(first considered by Ritt) of permutations with given cycle
structures. This idea leads to \fullref{even-deg_exceptions:thm} and
\fullref{mixed:deg:excep:cor}. It also allows us to establish the only
realizability result of this paper, \fullref{divisible:exist:thm}.
\end{itemize}

\paragraph{Comments on the new results} A very efficient algorithm
to treat the existence problem was developed in~\cite{Zheng} by
Zheng, who also produced a vast collection of experimental data
listing all the exceptional data up to degree $22$. The number of
these exceptions is very large, and it appears rather hard to
detect any sensible pattern in the list. However, we note that
over a half of these exceptions fall in the domain of
\fullref{even-deg_exceptions:thm} and a smaller, but still
noticeable, percentage is covered by
\fullref{excep:by:fixpoints:thm}.

We also notice that these two theorems cover all the exceptional
data with $d\leqslant 17$ and $n\geqslant 5$. In these cases, the
exceptions with even $d$ are almost always covered by
\fullref{even-deg_exceptions:thm}, although a few are
explained by \fullref{excep:by:fixpoints:thm}. At the level of
$d=18$ there are only two exceptional data with $n\geqslant 5$
which are not covered by these results, and a few more appear with
$19\leqslant d\leqslant 22$. To appreciate the power of these
statements, in particular for large $n$, notice also that the
total number of exceptional data in these degrees, even of those
with $n\geqslant 5$, is in the hundreds.

\paragraph{Comments on the new techniques}
A detailed account on the established methods for facing the
question of realizability of branch data will be given in
\fullref{review:section}, but we would like to mention here
that all these methods have a chiefly algebraic flavour, except
Bar\'anski's recent one.

The use we make of dessins d'enfants to prove
Propositions~\ref{non-realiz:53:prop} and~\ref{23:nonex:prop}, and
\fullref{excep:by:fixpoints:thm} provides the first
application of such a notion to the Hurwitz existence problem.
Besides the fact that the results are valuable on their own, we
consider it rather interesting to have a transparent
\emph{geometric} explanation of non-existence for infinite series
of coverings. In addition, the method appears to be generalizable
to more infinite series.

Even if very simple, our idea of analyzing coverings which should
be, if existent, compositions of other coverings proved rather
fruitful and also allowed us to explain non-realizability of many
data in purely geometric terms, as already stated above.

\paragraph{Organization of the paper}
In \fullref{review:section} we quickly review the results
previously known on the Hurwitz existence problem, with the aim of
helping the reader put our new contributions in the right
perspective. In \fullref{views:section} we describe the main
approaches historically taken to face the problem (and used in the
rest of our paper), also trying to make the relations between them
completely transparent. In \fullref{nonex:dessin:section} we
develop the technique of dessins d'enfants, which allows us to
establish Propositions~\ref{non-realiz:53:prop}
and~\ref{23:nonex:prop} and \fullref{excep:by:fixpoints:thm}.
In \fullref{nonex:decom:section} we apply the technique of
decomposability of coverings to prove
Theorems~\ref{even-deg_exceptions:thm}
and~\ref{divisible:exist:thm} and
\fullref{mixed:deg:excep:cor}.

\medskip
\noindent\textbf{Acknowledgements}\qua We are grateful to
Laurent~Bartholdi for his instructions on how to use the software
GAP to run computer experiments. We also thank Sergei Matveev,
Laura Mazzoni, Alexander Mednykh, and Joan Porti for helpful
conversations, and the referee of the first version of this paper
for very useful suggestions.

The first-named author was supported by the INTAS YS fellowship
03-55-1423. The second-named author was supported by the INTAS
project ``CalcoMet-GT'' 03-51-3663.

\section{Known results and techniques}\label{review:section}
In this section we outline the main partial solutions
of the Hurwitz existence problem which have been obtained
over the time.

\paragraph{Known results for $\Sigma\neq\matS$}
We begin by reviewing several results
whose overall content is that, to get a complete solution of the
Hurwitz existence problem, it would be sufficient to settle the
case where the base surface $\Sigma$ is the sphere $\matS$. The
first theorem we cite is attributed to Shephardson in Ezell~\cite[page
125]{Ezell} and explicitly proved in Husemoller~\cite[Theorem 4]{Husemoller}
and Edmonds, Kulkarni and Stong~\cite[Proposition 3.3]{EKS}:

\begin{thm}\label{OO:thm}
A compatible branch datum with $\Sigma$ orientable and
$\chi(\Sigma)\leqslant 0$ is realizable.
\end{thm}

The next result is proved in~\cite[Theorem 3.4]{Ezell}
and~\cite[Proposition 3.3]{EKS}:

\begin{thm}\label{NN:thm}
A compatible branch datum with $\Sigma$ and $\Sigmatil$
non-orientable and $\chi(\Sigma)\leqslant 0$ is realizable.
\end{thm}

We then quote the following elementary fact, stated
in~\cite[Proposition 2.7]{EKS}, and its
consequence~\cite[Proposition 3.4]{EKS}:

\begin{prop}\label{ON:reduction:prop}
A compatible branch datum with $\Sigma$ non-orientable and
$\Sigmatil$ orientable is realizable if and only if it is possible
to decompose for all $i$ the partition $(d_{ij})_{j=1,\ldots,m_i}$
of $d$ into partitions $(d'_{ij})_{j=1}^{m'_i}$ and
$(d''_{ij})_{j=1}^{m''_i}$ of $d/2$ in such a way that the branch
datum
$$\big(\Sigmatil,\Sigma',2n,d/2,(d'_{ij}),(d''_{ij})\big)$$
is realizable, where $\Sigma'$ is the orientable double covering
of $\Sigma$.
\end{prop}

\begin{cor}\label{ON:cor}
A compatible branch datum with $\Sigma$ non-orientable,
$\Sigmatil$ orientable and $\chi(\Sigma)\leqslant 0$ is
realizable.
\end{cor}

The next result is due to Edmonds, Kulkarni and
Stong~\cite[Theorem 5.1]{EKS}. We recall that $\matP$ is the
projective plane.

\begin{thm}\label{NP:thm}
A compatible branch datum with $\Sigma=\matP$ and non-orientable
$\Sigmatil$ is realizable.
\end{thm}

These theorems imply that only the following instances of the
Hurwitz existence problem remain open:
\begin{itemize} \item $\Sigma=\matS$;
\item $\Sigma=\matP$ and $\Sigmatil$ orientable.
\end{itemize}
However, \fullref{ON:reduction:prop} reduces the latter
instance to the former one, to which we will therefore confine
ourselves henceforth.

\paragraph{Known results for $\Sigma=\matS$}
When the base surface is the sphere $\matS$, not every compatible
branch datum is realizable.  The easiest example is given in degree
$d=4$ with $n=3$ branching points by the partitions
$(3,1),(2,2),(2,2)$, which implies that $\Sigmatil$ is $\matS$ too. A
proof follows from \fullref{EKS:222:prop} below, established in
\fullref{nonex:dessin:section}. In the rest of our paper the base
surface $\Sigma$ will always be the sphere $\matS$.  Recall that a
branch datum is \emph{exceptional} if it is compatible but not
realizable.

We will now review the main existence and non-existence theorems
proved in the literature.  However, we will not attempt to give a
comprehensive list of the abstract statements. Instead, we will
concentrate on the results which can be applied in a more direct
fashion. The following was established in~\cite[Proposition
5.7]{EKS}.

\begin{thm}\label{EKS:nonprime:thm}
For all non-prime $d$ there exist exceptional branch data of
degree $d$ with $n=3$ and $\Sigmatil=\matS$.
\end{thm}

Turning to existence, the most general known result appears to be
the following:

\begin{thm}\label{full:cycle:thm}
A compatible branch datum is realizable if one of the partitions
$(d_{ij})$ is given by $(d)$ only.
\end{thm}

This fact was first stated by Thom~\cite[Theorem 1]{Thom} for
$\Sigmatil$ also equal to $\matS$, was reproved by Khovanskii and
Zdravkovska~\cite[Theorem 2]{KZ} and Baranski~\cite[Theorem 6]{Baranski} in
the same context, and was generalized to arbitrary $\Sigmatil$
in~\cite[Proposition 5.2]{EKS}. A variation on this result is given
in~\cite[Proposition 5.3]{EKS}, where the realizable branch data with
one partition of the form $(d-1,1)$ are classified.  In~\cite[page
775]{EKS} a similar classification is announced for data with $n=3$
and one partition of the form $(m,1,\ldots,1)$.

Theorems~\ref{EKS:nonprime:thm} and~\ref{full:cycle:thm},
together with \fullref{EKS:222:prop} below,
are the main known results relevant to the case $n=3$, which we are most
interested in.
However, there are also some results relevant to the case where
$n$ is ``large'' (usually compared to $d$). In this respect, a
major contribution of Edmonds, Kulkarni and Stong is the
following~\cite[Theorem 5.4]{EKS}:

\begin{thm}\label{EKS:finiteness:thm} A branch datum with
$d\neq 4$ and $n\cdot d-\ntil\geqslant 3(d-1)$ is realizable. The
exceptional data with $d=4$ are precisely those with partitions
$(2,2),\ldots,(2,2),(3,1)$.
\end{thm}

Since we do not consider partitions of $d$ of the form
$(1,\ldots,1)$, one easily sees that $n\cdot d-\ntil\geqslant n$,
so a consequence of this result is that for fixed $d\neq 4$ the
number of exceptional branch data of degree $d$ is
finite~\cite[Corollary 4.4]{EKS}. The next result is established
in~\cite[Corollary 6.4]{EKS}:

\begin{prop}\label{EKS:222:prop}
For $\Sigmatil=\matS$, $n=3$ and even $d$ a branch datum with
partitions $(x,d-x),(2,\ldots,2),(2,\ldots,2)$ is realizable if
and only if $x=d/2$.
\end{prop}

The following results are due to Bar\'anski~\cite[Proposition
10, Theorem 12, Corollary 15]{Baranski}. We note that the first one
extends \fullref{full:cycle:thm} in the special case where
$\Sigmatil=\matS$, while the second one implies that for fixed $d$
the number of exceptional branch data of degree $d$ with
$\Sigmatil=\matS$ is finite (which was already known
after~\cite{EKS} for arbitrary $\Sigmatil$ and $d\neq 4$).

\begin{prop}\label{Baranski:long:cycles:prop}
A branch datum with $\Sigmatil=\matS$ is realizable if there
exists $r$ such that $m_1+\ldots+m_r=(r-1)\cdot d+1$.
\end{prop}

\begin{prop}\label{Baranski:finiteness:prop}
A compatible branch datum with $\Sigmatil=\matS$ and $n\geqslant
d$ is realizable.
\end{prop}

\begin{prop}\label{Baranski:small:prop}
A compatible branch datum with $\Sigmatil=\matS$, $d_{ij}\leqslant
2$ for all $i,j$, and $m_i\geqslant d-\sqrt{d/2}$ for all $i$ is
realizable.
\end{prop}

Considering \fullref{EKS:nonprime:thm} and the fact that all
known exceptional branch data occur with non-prime degree $d$, one
is naturally led to conjecture that for prime $d$ the Hurwitz
existence problem always has a positive solution. It is claimed
in~\cite[page 787]{EKS} that establishing this conjecture in the
special case $n=3$ would imply the general case.

\paragraph{Established techniques and related known results}
The proofs of Theorems~\ref{OO:thm} and~\ref{NN:thm} are based on
a reformulation of the main problem in terms of representations of
the fundamental group of the $n$--punctured surface $\Sigma$ into
the symmetric group $\permu_d$. This technique, already alluded to
several times above, goes back to Hurwitz
himself~\cite[Section~I.1]{Hurwitz} and was later revisited and
refined by Ezell~\cite[Theorem 2.1]{Ezell},
Singerman~\cite[Theorem 1]{Singerman}, where a more complicated
situation (which includes cusped surfaces) is considered, and
Husemoller~\cite[Theorems 3 and 5]{Husemoller}, where surfaces
with boundary are also accepted. The permutation technique will be
carefully reviewed and reinterpreted below in
\fullref{views:section}. It is however worth remarking here
that the main algebraic results on permutations leading to
Theorems~\ref{OO:thm} and~\ref{NN:thm} were also established in an
abstract context in~\cite{Boccara}.

In the special case where $\Sigma=\matS$ and $n=3$, the Hurwitz
existence problem can be reinterpreted in terms of algebraic maps
between algebraic curves. In this context Belyi's
theorem~\cite[Proposition 3]{Wolfart} could be viewed, in a sense,
as a solution of the problem, but the necessary and sufficient
condition for existence it gives is an abstract algebraic one
which it does not seem to be possible to check in practice.

The results of Bar\'anski stated above are based on a geometric
criterion~\cite[Lemma 5]{Baranski} for the realizability of a
branched covering of $\matS$ with $\Sigmatil$ also equal to
$\matS$.

To conclude, we quickly discuss the question (already mentioned
above) of counting the number of (suitably defined) equivalence
classes of branched coverings realizing a given branch datum. This
counting problem in fact goes back to Hurwitz as well, for whom it
was a primary motivation and who solved it
in~\cite[pages~16,~22]{Hurwitz} for the case $n=2d-2$ and all
partitions of the form $(2,1,\ldots,1)$, with $\Sigmatil=\matS$.
Complete formulae for the general case were given by
Mednykh~\cite[page 138]{Mednykh1},~\cite[Theorem C]{Mednykh2}. Of
course, a branch datum is realizable if and only if the
corresponding counting formula returns a positive value, so these
formulae give an implicit solution to the Hurwitz existence
problem. But again the actual computation appears to be hopeless.
A variation of Mednykh's formulae, where some explicit computation
is possible, was recently derived in~\cite{MSS} for special types
of branched coverings.

\paragraph{Computational aspects}
The realizability of any specific branch datum can in principle be
checked either using the formulae of Mednykh just mentioned, or
via Hurwitz's reformulation in terms of permutations. However the
former approach is too complicated to be practical. The latter one
is suited for computer implementation, but it involves operating
with conjugacy classes of permutations which already in the case
$n=3$ and $d=11$ are too huge for the capacity of today's
computers. A much more efficient method was recently developed by
Zheng in~\cite{Zheng}, where the existence problem was expressed
in terms of coefficients of certain generating functions. As
already mentioned, the formulae also established by Zheng for
calculating these coefficients allowed him to treat exceptional
data up to degree $22$, see~\cite{Zheng,Zheng-bis}.
\vspace{-2pt}

\section{Alternative viewpoints on the problem}\label{views:section}
\vspace{-2pt}

We start by describing equivalent formulations of the Hurwitz
existence problem in terms of ordinary coverings of surfaces with
boundary and in terms of 2--orbifolds. In some cases, one of these
reformulations will be easier to handle than the original
viewpoint of branched coverings.
\vspace{-2pt}

\paragraph{Surfaces with boundary}
Let $\Sigma_n$ denote the surface obtained from $\Sigma$ by
removing $n$ open discs with disjoint closures. A branched
covering realizing a datum
$\big(\Sigmatil,\Sigma,n,d,(d_{ij})_{i=1,\ldots,n}^{j=1,\ldots,m_i}\big)$
induces a genuine covering $\Sigmatil_{\ntil}\to\Sigma_n$ in which
the $i$-th component of $\partial\Sigma_n$ is covered by $m_i$
components of $\partial\Sigmatil_{\ntil}$, and the degrees of the
restrictions to these components are $(d_{ij})_{j=1,\ldots,m_i}$.
Conversely, any such genuine covering
$\Sigmatil_{\ntil}\to\Sigma_n$ induces a realization of the branch
datum. We will therefore consider also such a covering a
realization of the datum.
\vspace{-2pt}

\paragraph{2--orbifolds}
Let $(\Sigma;p_1,\ldots,p_n)$ denote the (closed locally
orientable) 2--orbifold with underlying surface $\Sigma$ and cone
points of orders $p_1,\ldots,p_n$. Choosing $p_i$ to be divisible
by $d_{ij}$ for $j=1,\ldots,m_i$, a branched covering realizing a
datum
$\big(\Sigmatil,\Sigma,n,d,(d_{ij})_{i=1,\ldots,n}^{j=1,\ldots,m_i}\big)$
induces an orbifold-covering~\cite{thurston:notes}
$$\left(\Sigmatil;\frac{p_1}{d_{11}},\ldots,\frac{p_1}{d_{1m_1}},\ldots,
\frac{p_n}{d_{n1}},\ldots,\frac{p_n}{d_{nm_n}}\right) \rightarrow
(\Sigma;p_1,\ldots,p_n).$$ And again, with details that we can
safely leave to the reader, a covering between 2--orbifolds induces
a branched covering between surfaces.
\vspace{-2pt}

\paragraph{Covering from permutations}
We will denote by $\matS$, $\matT$ and $\matP$ the sphere, the
torus and the projective plane respectively, whence by
$g\matT=\matT\#\ldots\#\matT$ ($g$ times) and
$g\matP=\matP\#\ldots\#\matP$ ($g$ times) the orientable and
non-orientable surfaces of genus $g\geqslant 1$. The following
group presentations are well-known:
\begin{eqnarray*}
\pi_1\big((g\matT)_n\big) & = &  \left\langle
a_1,b_1,\ldots,a_g,b_g,e_1,\ldots,e_n\big|\
[a_1,b_1]\cdots[a_g,b_g]\cdot e_1\cdots e_n \right\rangle\\
\pi_1\big((g\matP)_n\big) & = & \langle
c_1,\ldots,c_g,e_1,\ldots,e_n\big|\ c_1^2\cdots c_g^2\cdot
e_1\cdots e_n \rangle.
\end{eqnarray*}
In both cases the $e_i$'s are represented by the boundary circles.
In the orientable case $a_k,b_k$ is a meridian-longitude pair on
the $k$-th copy of $\matT$. In the non-orientable case $c_k$ is
the only non-trivial loop on the $k$-th copy of $\matP$.

The first alternative viewpoint on the Hurwitz existence problem
described above allows one to establish the following fact, originally
due to Hurwitz, Husemoller, Ezell, and Singerman. We quickly
review the geometric argument underlying the proof because we will
need it below.

\begin{thm}\label{Hurwitz:method:thm}
A compatible branch datum
$\big(\Sigmatil,\Sigma,n,d,(d_{ij})\big)$ is realizable if and
only if there exists a representation
$\theta\co \pi_1(\Sigma_n)\rightarrow \permu_d$ such that:
\begin{enumerate}
\item $\mbox{Im}(\theta)$ acts transitively on $\{1,\ldots,d\}$;
\item $\theta(e_i)$ has precisely $m_i$ cycles of lengths
$d_{i1},\ldots,d_{im_i}$; \item For non-orientable $\Sigma$ and
orientable $\Sigmatil$, each permutation $\theta(c_k)$ consists of
cycles of even length only.
\end{enumerate}
\end{thm}

\begin{proof}
Our proof only works for $n>0$, ie, it does not for
genuine coverings. Recall first~\cite[Chapter 2]{MaFo} that
$\Sigma_n$ can be obtained as shown in
\fullref{discs-with-gluings:fig}
    \begin{figure}[ht!]\hspace{-9mm}
\begin{picture}(0,0)%
\includegraphics[scale=1.11]{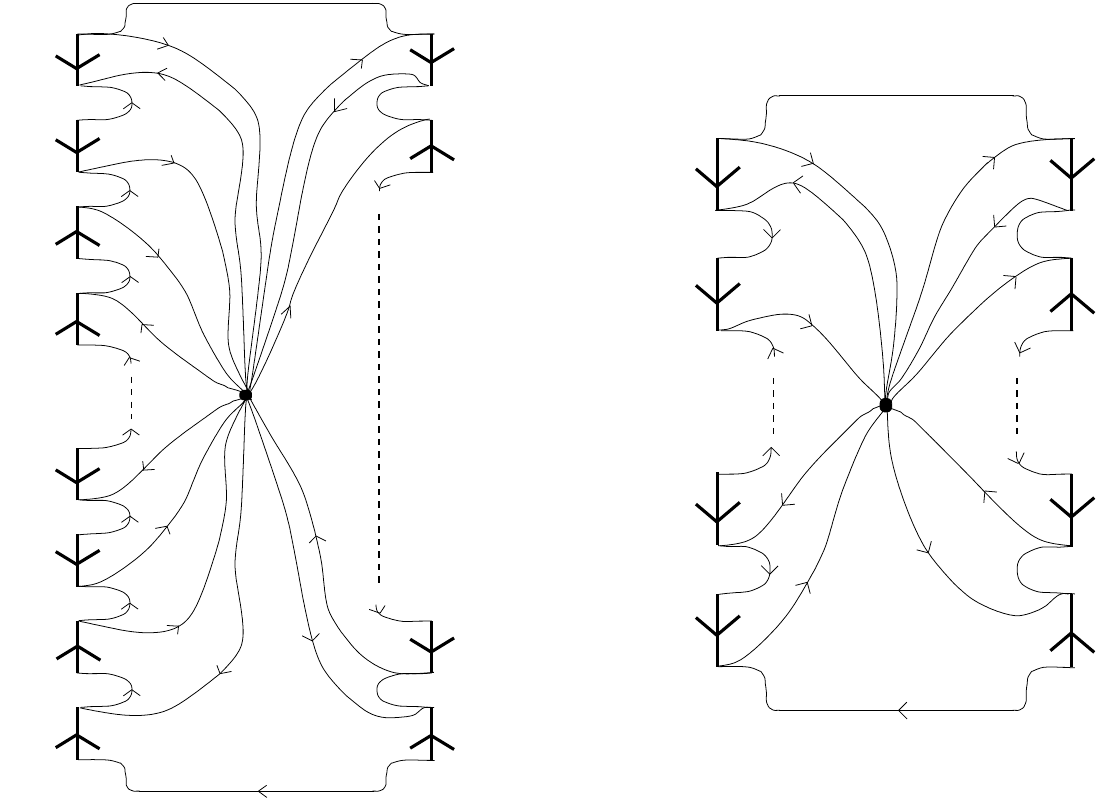}%
\end{picture}%
\setlength{\unitlength}{2190sp}%
\def\SetFigFont#1#2#3#4#5{\footnotesize}%
\begin{picture}(10538,7649)(601,-6662)
\put(2251,614){\makebox(0,0)[lb]{\smash{\SetFigFont{6}{7.2}{\familydefault}{\mddefault}{\updefault}{$e_n$}%
}}}
\put(4351,-886){\makebox(0,0)[lb]{\smash{\SetFigFont{6}{7.2}{\familydefault}{\mddefault}{\updefault}{$e_n$}%
}}}
\put(2176,-1411){\makebox(0,0)[lb]{\smash{\SetFigFont{6}{7.2}{\familydefault}{\mddefault}{\updefault}{$a_1$}%
}}}
\put(601,389){\makebox(0,0)[lb]{\smash{\SetFigFont{6}{7.2}{\familydefault}{\mddefault}{\updefault}{$\alpha_{1,+}$}%
}}}
\put(601,-436){\makebox(0,0)[lb]{\smash{\SetFigFont{6}{7.2}{\familydefault}{\mddefault}{\updefault}{$\beta_{1,+}$}%
}}}
\put(601,-2086){\makebox(0,0)[lb]{\smash{\SetFigFont{6}{7.2}{\familydefault}{\mddefault}{\updefault}{$\beta_{1,-}$}%
}}}
\put(1501,-811){\makebox(0,0)[lb]{\smash{\SetFigFont{6}{7.2}{\familydefault}{\mddefault}{\updefault}{$e_n$}%
}}}
\put(1501, 14){\makebox(0,0)[lb]{\smash{\SetFigFont{6}{7.2}{\familydefault}{\mddefault}{\updefault}{$e_n$}%
}}}
\put(1501,-1636){\makebox(0,0)[lb]{\smash{\SetFigFont{6}{7.2}{\familydefault}{\mddefault}{\updefault}{$e_n$}%
}}}
\put(1501,-2461){\makebox(0,0)[lb]{\smash{\SetFigFont{6}{7.2}{\familydefault}{\mddefault}{\updefault}{$e_n$}%
}}}
\put(3676,464){\makebox(0,0)[lb]{\smash{\SetFigFont{6}{7.2}{\familydefault}{\mddefault}{\updefault}{$e_n$}%
}}}
\put(4951,389){\makebox(0,0)[lb]{\smash{\SetFigFont{6}{7.2}{\familydefault}{\mddefault}{\updefault}{$\epsilon_{n-1,+}$}%
}}}
\put(4951,-436){\makebox(0,0)[lb]{\smash{\SetFigFont{6}{7.2}{\familydefault}{\mddefault}{\updefault}{$\epsilon_{n-1,-}$}%
}}}
\put(2101, 89){\makebox(0,0)[lb]{\smash{\SetFigFont{6}{7.2}{\familydefault}{\mddefault}{\updefault}{$a_1$}%
}}}
\put(2176,-436){\makebox(0,0)[lb]{\smash{\SetFigFont{6}{7.2}{\familydefault}{\mddefault}{\updefault}{$b_1$}%
}}}
\put(1951,-2011){\makebox(0,0)[lb]{\smash{\SetFigFont{6}{7.2}{\familydefault}{\mddefault}{\updefault}{$b_1$}%
}}}
\put(601,-3586){\makebox(0,0)[lb]{\smash{\SetFigFont{6}{7.2}{\familydefault}{\mddefault}{\updefault}{$\alpha_{g,+}$}%
}}}
\put(601,-4411){\makebox(0,0)[lb]{\smash{\SetFigFont{6}{7.2}{\familydefault}{\mddefault}{\updefault}{$\beta_{g,+}$}%
}}}
\put(601,-6211){\makebox(0,0)[lb]{\smash{\SetFigFont{6}{7.2}{\familydefault}{\mddefault}{\updefault}{$\beta_{g,-}$}%
}}}
\put(3376,-5311){\makebox(0,0)[lb]{\smash{\SetFigFont{6}{7.2}{\familydefault}{\mddefault}{\updefault}{$e_1$}%
}}}
\put(3751,-4111){\makebox(0,0)[lb]{\smash{\SetFigFont{6}{7.2}{\familydefault}{\mddefault}{\updefault}{$e_1$}%
}}}
\put(1951,-3736){\makebox(0,0)[lb]{\smash{\SetFigFont{6}{7.2}{\familydefault}{\mddefault}{\updefault}{$a_g$}%
}}}
\put(2176,-4336){\makebox(0,0)[lb]{\smash{\SetFigFont{6}{7.2}{\familydefault}{\mddefault}{\updefault}{$b_g$}%
}}}
\put(2701,-5686){\makebox(0,0)[lb]{\smash{\SetFigFont{6}{7.2}{\familydefault}{\mddefault}{\updefault}{$b_g$}%
}}}
\put(5026,-5311){\makebox(0,0)[lb]{\smash{\SetFigFont{6}{7.2}{\familydefault}{\mddefault}{\updefault}{$\epsilon_{1,+}$}%
}}}
\put(5026,-6136){\makebox(0,0)[lb]{\smash{\SetFigFont{6}{7.2}{\familydefault}{\mddefault}{\updefault}{$\epsilon_{1,-}$}%
}}}
\put(4351,-4861){\makebox(0,0)[lb]{\smash{\SetFigFont{6}{7.2}{\familydefault}{\mddefault}{\updefault}{$e_n$}%
}}}
\put(3440,-2124){\makebox(0,0)[lb]{\smash{\SetFigFont{6}{7.2}{\familydefault}{\mddefault}{\updefault}{$e_{n-1}$}%
}}}
\put(3774,-234){\makebox(0,0)[lb]{\smash{\SetFigFont{6}{7.2}{\familydefault}{\mddefault}{\updefault}{$e_{n-1}$}%
}}}
\put(1501,-3211){\makebox(0,0)[lb]{\smash{\SetFigFont{6}{7.2}{\familydefault}{\mddefault}{\updefault}{$e_n$}%
}}}
\put(1501,-3961){\makebox(0,0)[lb]{\smash{\SetFigFont{6}{7.2}{\familydefault}{\mddefault}{\updefault}{$e_n$}%
}}}
\put(1576,-4786){\makebox(0,0)[lb]{\smash{\SetFigFont{6}{7.2}{\familydefault}{\mddefault}{\updefault}{$e_n$}%
}}}
\put(1576,-5611){\makebox(0,0)[lb]{\smash{\SetFigFont{6}{7.2}{\familydefault}{\mddefault}{\updefault}{$e_n$}%
}}}
\put(2251,-5236){\makebox(0,0)[lb]{\smash{\SetFigFont{6}{7.2}{\familydefault}{\mddefault}{\updefault}{$a_g$}%
}}}
\put(7576,-4486){\makebox(0,0)[lb]{\smash{\SetFigFont{6}{7.2}{\familydefault}{\mddefault}{\updefault}{$e_n$}%
}}}
\put(7651,-2386){\makebox(0,0)[lb]{\smash{\SetFigFont{6}{7.2}{\familydefault}{\mddefault}{\updefault}{$e_n$}%
}}}
\put(8476,-511){\makebox(0,0)[lb]{\smash{\SetFigFont{6}{7.2}{\familydefault}{\mddefault}{\updefault}{$e_n$}%
}}}
\put(10501,-3436){\makebox(0,0)[lb]{\smash{\SetFigFont{6}{7.2}{\familydefault}{\mddefault}{\updefault}{$e_n$}%
}}}
\put(10501,-2461){\makebox(0,0)[lb]{\smash{\SetFigFont{6}{7.2}{\familydefault}{\mddefault}{\updefault}{$e_n$}%
}}}
\put(9751,-436){\makebox(0,0)[lb]{\smash{\SetFigFont{6}{7.2}{\familydefault}{\mddefault}{\updefault}{$e_n$}%
}}}
\put(8251,-1036){\makebox(0,0)[lb]{\smash{\SetFigFont{6}{7.2}{\familydefault}{\mddefault}{\updefault}{$c_1$}%
}}}
\put(8176,-4036){\makebox(0,0)[lb]{\smash{\SetFigFont{6}{7.2}{\familydefault}{\mddefault}{\updefault}{$c_g$}%
}}}
\put(8401,-4786){\makebox(0,0)[lb]{\smash{\SetFigFont{6}{7.2}{\familydefault}{\mddefault}{\updefault}{$c_g$}%
}}}
\put(9301,-4486){\makebox(0,0)[lb]{\smash{\SetFigFont{6}{7.2}{\familydefault}{\mddefault}{\updefault}{$e_1$}%
}}}
\put(9751,-3886){\makebox(0,0)[lb]{\smash{\SetFigFont{6}{7.2}{\familydefault}{\mddefault}{\updefault}{$e_1$}%
}}}
\put(6751,-811){\makebox(0,0)[lb]{\smash{\SetFigFont{6}{7.2}{\familydefault}{\mddefault}{\updefault}{$\gamma_{1,+}$}%
}}}
\put(6751,-1936){\makebox(0,0)[lb]{\smash{\SetFigFont{6}{7.2}{\familydefault}{\mddefault}{\updefault}{$\gamma_{1,-}$}%
}}}
\put(6751,-4036){\makebox(0,0)[lb]{\smash{\SetFigFont{6}{7.2}{\familydefault}{\mddefault}{\updefault}{$\gamma_{g,+}$}%
}}}
\put(6751,-5236){\makebox(0,0)[lb]{\smash{\SetFigFont{6}{7.2}{\familydefault}{\mddefault}{\updefault}{$\gamma_{g,-}$}%
}}}
\put(11101,-736){\makebox(0,0)[lb]{\smash{\SetFigFont{6}{7.2}{\familydefault}{\mddefault}{\updefault}{$\epsilon_{n-1,+}$}%
}}}
\put(11101,-1861){\makebox(0,0)[lb]{\smash{\SetFigFont{6}{7.2}{\familydefault}{\mddefault}{\updefault}{$\epsilon_{n-1,-}$}%
}}}
\put(11101,-4036){\makebox(0,0)[lb]{\smash{\SetFigFont{6}{7.2}{\familydefault}{\mddefault}{\updefault}{$\epsilon_{1,+}$}%
}}}
\put(11101,-5161){\makebox(0,0)[lb]{\smash{\SetFigFont{6}{7.2}{\familydefault}{\mddefault}{\updefault}{$\epsilon_{1,-}$}%
}}}
\put(10081,-2011){\makebox(0,0)[lb]{\smash{\SetFigFont{6}{7.2}{\familydefault}{\mddefault}{\updefault}{$e_{n-1}$}%
}}}
\put(7651,-3361){\makebox(0,0)[lb]{\smash{\SetFigFont{6}{7.2}{\familydefault}{\mddefault}{\updefault}{$e_n$}%
}}}
\put(3076,-6436){\makebox(0,0)[lb]{\smash{\SetFigFont{6}{7.2}{\familydefault}{\mddefault}{\updefault}{$e_n$}%
}}}
\put(601,-5236){\makebox(0,0)[lb]{\smash{\SetFigFont{6}{7.2}{\familydefault}{\mddefault}{\updefault}{$\alpha_{g,-}$}%
}}}
\put(601,-1261){\makebox(0,0)[lb]{\smash{\SetFigFont{6}{7.2}{\familydefault}{\mddefault}{\updefault}{$\alpha_{1,-}$}%
}}}
\put(8476,-2086){\makebox(0,0)[lb]{\smash{\SetFigFont{6}{7.2}{\familydefault}{\mddefault}{\updefault}{$c_1$}%
}}}
\put(9076,-5611){\makebox(0,0)[lb]{\smash{\SetFigFont{6}{7.2}{\familydefault}{\mddefault}{\updefault}{$e_n$}%
}}}
\put(9916,-834){\makebox(0,0)[lb]{\smash{\SetFigFont{6}{7.2}{\familydefault}{\mddefault}{\updefault}{$e_{n-1}$}%
}}}
\put(7651,-1261){\makebox(0,0)[lb]{\smash{\SetFigFont{6}{7.2}{\familydefault}{\mddefault}{\updefault}{$e_n$}%
}}}
\end{picture}
\caption{Construction of punctured surfaces} \label{discs-with-gluings:fig}
\vspace{10pt}
\end{figure}
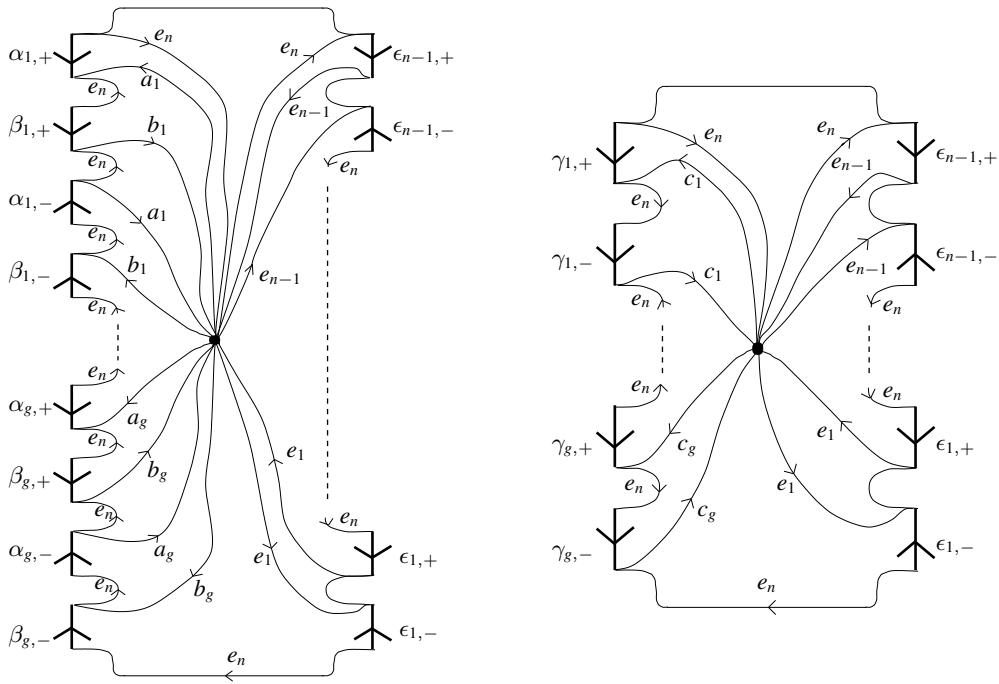
from a disc $\Delta$ by identifying each arc $\epsilon_{i,-}$ to
the arc $\epsilon_{i,+}$ and, depending on orientability, either
each $\alpha_{k,-}$ to $\alpha_{k,+}$ and $\beta_{k,-}$ to
$\beta_{k,+}$, or each $\gamma_{k,-}$ to $\gamma_{k,+}$. The
figure also shows loops corresponding to the generators of
$\pi_1(\Sigma_n)$ used in the above presentations.

The degree-$d$ coverings of $\Sigma_n$ are now obtained as
follows. We first take the disjoint union of $d$ copies
$(\Delta^{(h)})_{h=1}^d$ of $\Delta$, with the corresponding arcs
$\epsilon^{(h)}_{i,\pm}$ etc. Then we glue each
$\epsilon^{(h)}_{i,-}$ to some
$\epsilon^{(\theta(e_i)(h))}_{i,+}$, where
$\theta(e_i)\in\permu_d$, and similarly for the other arcs, using
permutations $\theta(a_k)$ and $\theta(b_k)$, or $\theta(c_k)$.
The corresponding covering is induced by the identification of
each $\Delta^{(h)}$ with $\Delta$, and of course the covering is
connected if and only if the subgroup of $\permu_d$ generated by
$\theta(e_1),\ldots,\theta(e_{n-1})$ and either the
$\theta(a_k)$'s and $\theta(b_k)$'s or the $\theta(c_k)$'s acts
transitively on $\{1,\ldots,d\}$. It is also obvious that for
$i=1,\ldots,n-1$ the way the $i$-th component of
$\partial\Sigma_n$ is covered depends on the cyclic structure of
$\theta(e_i)$ as described in the statement.

Considering the form of the presentation of $\pi_1(\Sigma_n)$
given above, we see that we can define $\theta(e_n)\in\permu_d$ in
a unique fashion so to get a representation $\theta$. To conclude
we must show that the way the $n$-th component of
$\partial\Sigma_n$ is covered depends on the cyclic structure of
$\theta(e_n)$, and that the covering over non-orientable
$\Sigma_n$ is orientable if and only if each $\theta(c_k)$
contains cycles of even length only. Both assertions are easy and
left as an exercise.
\end{proof}

\paragraph{The parity condition}
Using \fullref{Hurwitz:method:thm}
we can now explain Condition~\ref{Pcond} of
the definition of compatible branch datum.

\begin{lemma}
If a datum is realizable then $n\cdot d-\ntil$ is even.
\end{lemma}

\begin{proof}
The conclusion is evident when $\Sigma$ is orientable, so we assume
it is not. Consider a representation $\theta$
realizing the datum as in \fullref{Hurwitz:method:thm}. Notice
that $\theta(e_1)\cdots\theta(e_n)$ is an even permutation, since
its inverse is the product of either commutators or squares. A
permutation is even if and only if it contains an even number of
cycles of even length, ie, if the sum of a contribution
$\ell-1$ for each cycle of length $\ell$ is even. Now the lengths
of the cycles of $\theta(e_i)$ are the $d_{ij}$'s, whence
$$\sum_{i=1}^{n}\sum_{j=1}^{m_i}(d_{ij}-1)=n\cdot d-\ntil$$
is even.
\end{proof}

\section{Exceptions via dessins d'enfants}\label{nonex:dessin:section}
As mentioned in \fullref{new:results:section} many exceptional
branch data exist when $\Sigma$ is the sphere $\matS$, and in this
section we present several classes of them, using a variation on 
Grothendieck's dessins d'enfants~\cite{Groth}.  This notion in
its original form is only relevant to the case of $n=3$ branching
points, but we actually generalize it to arbitrary $n\geqslant 3$.

\begin{defn}
\emph{A \emph{dessin d'enfant} on
$\Sigmatil$ is a graph $D\subset\Sigmatil$ where:
\begin{enumerate}
\item For some $n\geqslant 3$ the set of vertices of $D$ is split
as $V_1\sqcup\ldots\sqcup V_{n-1}$ and the set of edges of $D$ is
split as $E_1\sqcup\ldots\sqcup E_{n-2}$; \item For
$i=1,\ldots,n-2$ each edge in $E_i$ joins a vertex of $V_i$ to one of
$V_{i+1}$; \item For $i=2,\ldots,n-2$ any vertex of $V_i$ has even
valence and going around it we alternatively encounter edges from
$E_{i-1}$ and edges from $E_{i}$; \item $\Sigmatil\setminus D$
consists of open discs.
\end{enumerate}
The \emph{length} of one of the discs in $\Sigmatil\setminus D$
is the number of edges of $D$ along which the boundary of the disc
passes (with multiplicity).}
\end{defn}

\begin{prop}\label{dessins:coverings:prop}
The realizations of a branch datum
$\big(\Sigmatil,\matS,n,d,(d_{ij})\big)$ correspond to the
dessins d'enfants $D\subset\Sigmatil$ with the set of
vertices split as $V_1\sqcup V_2\sqcup\ldots\sqcup V_{n-1}$ such
that for $i=1$ and $i=n-1$ the vertices in $V_i$ have valences
$(d_{ij})_{j=1,\ldots,m_i}$, for $i=2,\ldots,n-2$ the vertices in
$V_i$ have valences $(2d_{ij})_{j=1,\ldots,m_i}$, and the discs in
$\Sigmatil\setminus D$ have lengths
$(2(n-2)d_{nj})_{j=1,\ldots,m_n}$.
\end{prop}

\begin{proof}
Suppose a realization $f\co \Sigmatil\to\matS$ exists, let the
branching points be $x_1,\ldots,x_n$, for $i=1,\ldots,n-2$ choose
a simple arc $\alpha_i$ joining $x_i$ to $x_{i+1}$, suppose the
$\alpha_i$'s meet at their ends only and avoid $x_n$, and let
$\alpha$ be their union. Then define $D$ as $f^{-1}(\alpha)$ and
set $V_i=f^{-1}(x_i)$ and $E_i=f^{-1}(\alpha_i)$. To conclude that
$D$ is a dessin d'enfant with valences and lengths as required,
the only non-obvious facts concern the components of
$\Sigmatil\setminus D$. But $\matS\setminus\alpha$ is an open
disc, the restriction of $f$ to any component of
$\Sigmatil\setminus D$ is a covering onto this disc with a single
branching point, and such a covering is always modelled on the
covering $z\mapsto z^k$ of the open unit disc onto itself, so the
components of $\Sigmatil\setminus D$ are open discs. More
precisely, there is one such disc for each element of
$f^{-1}(x_n)$, and it is easy to see that the $j$-th one has
length as required.

Reversing this construction is a routine matter left to the
reader.\end{proof}

\paragraph{From dessins to permutations and back}
We recall that the conjugacy classes in $\permu_d$ are given
precisely by the partitions of $d$, with the class of a
permutation being the array of lengths of its cycles. The
following is a consequence of \fullref{Hurwitz:method:thm}, in
which we emphasize the constructive nature of the theorem:

\begin{cor}\label{permu:S:n:cor}
The realizations of a branch datum
$\big(\Sigmatil,\matS,n,d,(d_{ij})\big)$ correspond to the choices
of $\tau_i=\theta(e_i)$ in $\permu_d$ for $i=1,\ldots,n-1$
such that $\langle \tau_1,\ldots,\tau_{n-1}\rangle$ is transitive and,
setting $\tau_n=\tau_{n-1}^{-1}\cdots\tau_1^{-1}$, for $i=1,\ldots,n$ the
conjugacy class of $\tau_i$ is given by
$(d_{ij})_{j=1,\ldots,m_i}$.
\end{cor}

Combining \fullref{dessins:coverings:prop} and
\fullref{permu:S:n:cor} we deduce a correspondence between
(suitable) dessins d'enfants in $\Sigmatil$ and (suitable) choices
of $\tau_1,\ldots,\tau_{n-1}\in\permu_d$. Since we will use it in the sequel,
we spell out this correspondence explicitly in the next two propositions.
Proofs are easy and hence omitted.

\begin{prop}\label{dessins:to:permus:prop}
Given a dessin d'enfant $D$, with notation as in the definition,
corresponding to a realization of a branch datum
$\big(\Sigmatil,\matS,n,d,(d_{ij})\big)$, permutations\break
$\tau_1,\ldots,\tau_{n-1}$ corresponding to the same realization
are constructed as follows:

\begin{itemize}
\item Enumerate the edges of $E_i$ as
$e_i^{(1)},\ldots,e_i^{(d)}$, starting in an arbitrary fashion for
$E_1$ and so that for $i\geqslant 2$ around each vertex of $V_{i}$
each edge $e_{i}^{(k)}$ is followed by the edge $e_{i-1}^{(k)}$
with the same number $k$;

\item For $i\leqslant n-2$ and $k\in\{1,\ldots,d\}$ select the
vertex of $V_i$ to which $e_i^{(k)}$ is incident and define
$\tau_i(k)$ to be $h$ such that the next $e_i^{(*)}$ around the
vertex is $e_i^{(h)}$.

\item For $k\in\{1,\ldots,d\}$ select the vertex of $V_{n-1}$ to
which $e_{n-2}^{(k)}$ is incident and define $\tau_{n-1}(k)$ to be
$h$ such that the next $e_{n-2}^{(*)}$ around the vertex is
$e_{n-2}^{(h)}$.
\end{itemize}
\end{prop}

Concerning this statement, recall that $\Sigmatil$ is oriented and
note the following about the edge $e_i^{(h)}$ needed to define
$\tau_i(k)$: for $i=1$ this edge comes immediately after
$e_i^{(k)}$, while for $i=2,\ldots,n-2$ there is the edge
$e_{i-1}^{(k)}$ located in between. Notice also that the edge
$e_{n-2}^{(h)}$ when used in the definition of $\tau_{n-1}(k)$
again comes immediately after $e_{n-2}^{(k)}$.

To describe the opposite correspondence, given
$\tau_1,\ldots,\tau_{n-1}\in\permu_d$, we construct a graph
$D(\tau_1,\ldots,\tau_{n-1})$ with vertices $V_1\sqcup\ldots\sqcup
V_{n-1}$, where $V_i$ is the set of cycles of $\tau_i$, and for
$i=1,\ldots,n-2$ and $k=1,\ldots,d$ an edge $e_i^{(k)}$ joins the
cycles of $\tau_i$ and $\tau_{i+1}$ which contain $k$. Note that
$D(\tau_1,\ldots,\tau_{n-1})$ is connected if and only if
$\langle\tau_1,\ldots,\tau_{n-1}\rangle$ is transitive. Then we
consider in $D(\tau_1,\ldots,\tau_{n-1})$ the loops constructed as
follows: we start with some vertex of $V_1$ and follow the path
$e_1^{(k)},\ldots,e_{n-2}^{(k)}$. Having thus arrived to a vertex
of $V_{n-1}$, we follow the path
$$e_{n-2}^{(\tau_{n-1}^{-1}(k))},
e_{n-3}^{(\tau_{n-2}^{-1}\tau_{n-1}^{-1}(k))},
\ldots,e_2^{(\tau_3^{-1}\cdots\tau_{n-1}^{-1}(k))},
e_1^{(\tau_2^{-1}\tau_3^{-1}\ldots\tau_{n-1}^{-1}(k))}$$ which
takes us to a vertex of $V_1$. From there we proceed similarly
starting from
$e_1^{(\tau_1^{-1}\tau_2^{-1}\tau_3^{-1}\ldots\tau_{n-1}^{-1}(k))}$
until we find the edge $e_1^{(k)}$ again.

\begin{prop}\label{permu:to:ext:dessin:prop}
Given $\tau_1,\ldots,\tau_{n-1}\in\permu_d$ corresponding to a
realization of a branch datum
$\big(\Sigmatil,\matS,n,d,(d_{ij})\big)$, the space obtained from
$D(\tau_1,\ldots,\tau_{n-1})$ by attaching discs to the loops just
described is $\Sigmatil$, and
$D(\tau_1,\ldots,\tau_{n-1})\subset\Sigmatil$ is a
dessin d'enfant corresponding to the same realization of the
branch datum.
\end{prop}

\paragraph{A sample application}
To investigate the realizability of a given branch datum\break
$\big(\Sigmatil,\matS,3,d,(d_{ij})\big)$ using the permutation
approach of \fullref{permu:S:n:cor}, one should fix a
certain $\tau_1$ with cycle lengths $(d_{1j})$ and then let
$\tau_2$ vary in the conjugacy class $(d_{2j})$, checking that
$\langle\tau_1,\tau_2\rangle$ is transitive and that
$\tau_1\cdot\tau_2$ has cycle lengths $(d_{3j})$. Since conjugacy
classes are huge, this method is only feasible for very small $d$,
and we exploited it using the software GAP for $d\leqslant 10$
(but note that Zheng's alternative method~\cite{Zheng} allows one
to treat much higher degrees). On the other hand, the
realizability criterion through dessins d'enfants, stated in
\fullref{dessins:coverings:prop}, has the advantage of
usually requiring the consideration of a much smaller number of
cases. The geometric nature of the criterion often also makes it
very easy to apply it. As a first example, we give a very simple
proof of a result stated in \fullref{review:section} and
originally established in~\cite{EKS}.

\dimo{EKS:222:prop} If the datum is realizable then the dessin
d'enfant associated to the last two
branching points is just a circle embedded in $\matS$. Such a
dessin decomposes $\matS$ into two discs of the same length, so
$x=d-x$, whence $x=d/2$. The same argument proves also the
opposite implication. \finedimo

\paragraph{Graph fattening and applications}
To apply \fullref{dessins:coverings:prop} it is sometimes
useful to switch the viewpoint: instead of trying to embed a
dessin $D$ in the surface $\Sigmatil$, we try to thicken a given
graph $D$ to a surface with boundary so to get $\Sigmatil$ by
capping off the boundary circles.  An application of this method
is given by the following proof of one of the results stated in
\fullref{new:results:section}.

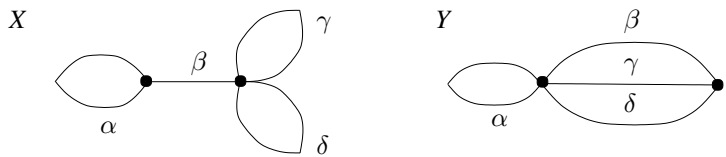
\begin{figure}[ht!]
\centering
\begin{picture}(0,0)%
\includegraphics[scale=0.9]{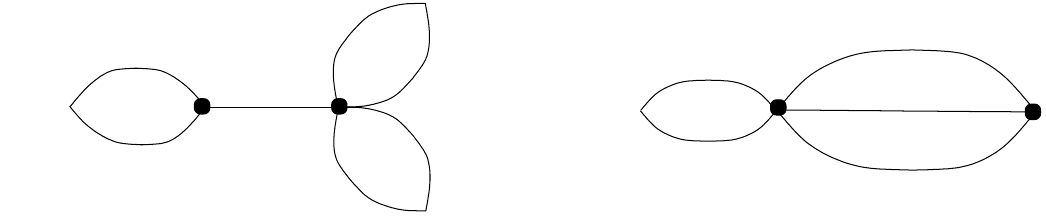}%
\end{picture}%
\setlength{\unitlength}{2308sp}%
\def\SetFigFont#1#2#3#4#5{\small}%
\begin{picture}(7697,1638)(226,-283)
\put(3526,-211){\makebox(0,0)[lb]{\smash{\SetFigFont{12}{14.4}{\familydefault}{\mddefault}{\updefault}{$\delta$}%
}}}
\put(226,1139){\makebox(0,0)[lb]{\smash{\SetFigFont{12}{14.4}{\familydefault}{\mddefault}{\updefault}{$X$}%
}}}
\put(3526,1139){\makebox(0,0)[lb]{\smash{\SetFigFont{12}{14.4}{\familydefault}{\mddefault}{\updefault}{$\gamma$}%
}}}
\put(2176,689){\makebox(0,0)[lb]{\smash{\SetFigFont{12}{14.4}{\familydefault}{\mddefault}{\updefault}{$\beta$}%
}}}
\put(1201, 14){\makebox(0,0)[lb]{\smash{\SetFigFont{12}{14.4}{\familydefault}{\mddefault}{\updefault}{$\alpha$}%
}}}
\put(5401, 89){\makebox(0,0)[lb]{\smash{\SetFigFont{12}{14.4}{\familydefault}{\mddefault}{\updefault}{$\alpha$}%
}}}
\put(6826,1139){\makebox(0,0)[lb]{\smash{\SetFigFont{12}{14.4}{\familydefault}{\mddefault}{\updefault}{$\beta$}%
}}}
\put(6826,689){\makebox(0,0)[lb]{\smash{\SetFigFont{12}{14.4}{\familydefault}{\mddefault}{\updefault}{$\gamma$}%
}}}
\put(6826,239){\makebox(0,0)[lb]{\smash{\SetFigFont{12}{14.4}{\familydefault}{\mddefault}{\updefault}{$\delta$}%
}}}
\put(4801,1139){\makebox(0,0)[lb]{\smash{\SetFigFont{12}{14.4}{\familydefault}{\mddefault}{\updefault}{$Y$}%
}}}
\end{picture}
    \caption{The two graphs with vertices of valences 5 and 3} \label{graphs5-3:fig}
    \end{figure}

\dimo{non-realiz:53:prop} We give a unified proof. A dessin
d'enfant $D$ corresponding to the first two partitions in the
given branch datum, as an abstract graph, is homeomorphic to one
of the graphs $X$ and $Y$ shown in \fullref{graphs5-3:fig},
with the two visible vertices lying in $V_2$. To get $D$ we must
then insert on each edge of $X$ or $Y$ an odd number of vertices
belonging to $V_1$ and $V_2$ alternatively. With a slight abuse of
notation, suppose we add $2\alpha+1$ vertices on $\alpha$, then
$2\beta+1$ on $\beta$, and so on. Using the fact that $V_1$ has
$d/2$ vertices, we see that $\alpha+\beta+\gamma+\delta=d/2-4$,
and this is the only constraint on $\alpha,\beta,\gamma,\delta$.

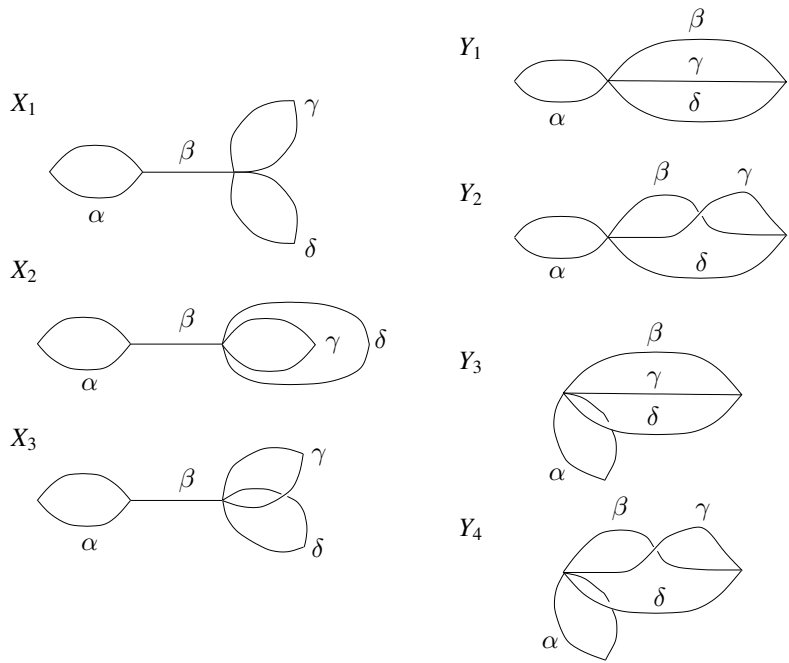
\begin{figure}[ht!]
    \begin{center}
\begin{picture}(0,0)%
\includegraphics[scale=0.9]{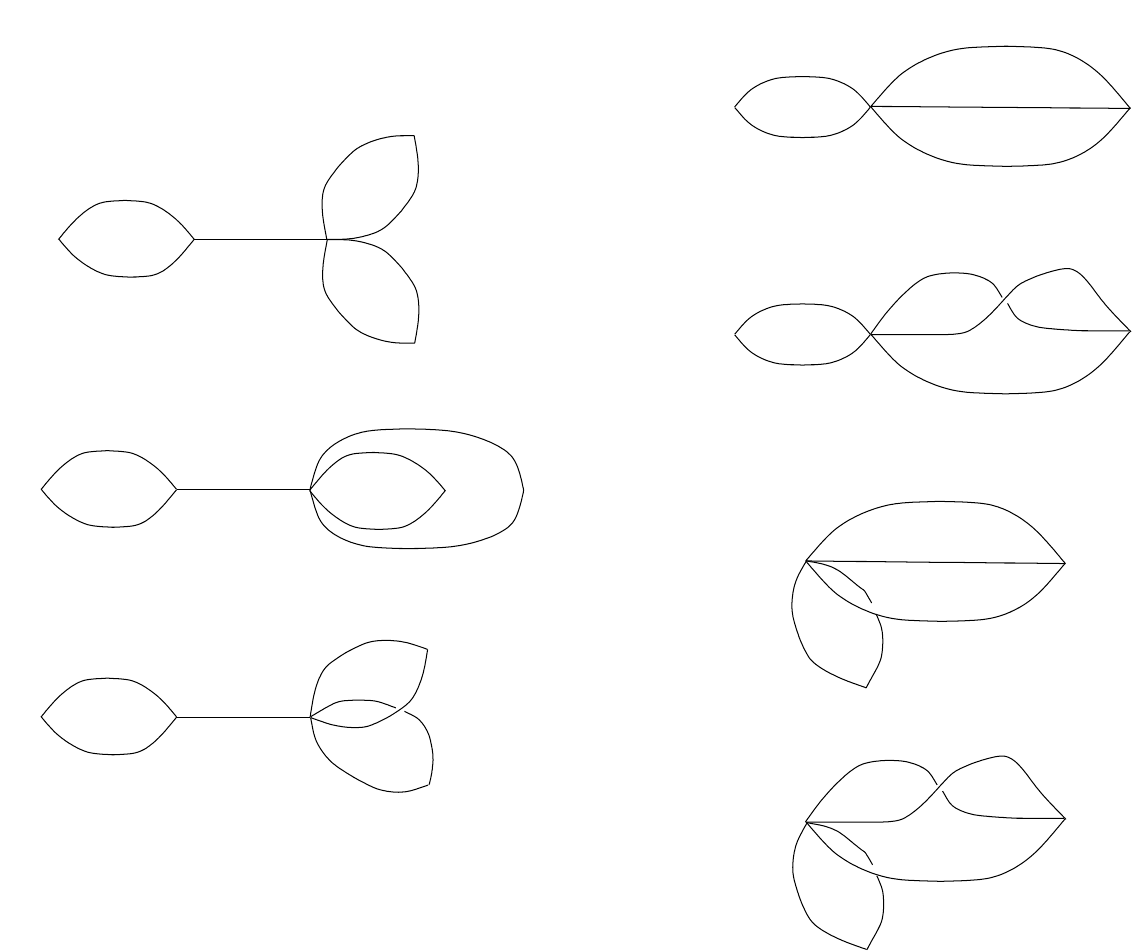}%
\end{picture}%
\setlength{\unitlength}{2308sp}%
\begingroup\makeatletter\ifx\SetFigFont\undefined%
\gdef\SetFigFont#1#2#3#4#5{%
  \reset@font\fontsize{#1}{#2pt}%
  \fontfamily{#3}\fontseries{#4}\fontshape{#5}%
  \selectfont}%
\fi\endgroup%
\begin{picture}(8364,7012)(601,-5306)
\put(3751,-961){\makebox(0,0)[lb]{\smash{\SetFigFont{10}{12.0}{\familydefault}{\mddefault}{\updefault}{$\delta$}%
}}}
\put(3826,-3136){\makebox(0,0)[lb]{\smash{\SetFigFont{10}{12.0}{\familydefault}{\mddefault}{\updefault}{$\gamma$}%
}}}
\put(3826,-4186){\makebox(0,0)[lb]{\smash{\SetFigFont{10}{12.0}{\familydefault}{\mddefault}{\updefault}{$\delta$}%
}}}
\put(2401, 89){\makebox(0,0)[lb]{\smash{\SetFigFont{10}{12.0}{\familydefault}{\mddefault}{\updefault}{$\beta$}%
}}}
\put(2401,-1711){\makebox(0,0)[lb]{\smash{\SetFigFont{10}{12.0}{\familydefault}{\mddefault}{\updefault}{$\beta$}%
}}}
\put(2401,-3436){\makebox(0,0)[lb]{\smash{\SetFigFont{10}{12.0}{\familydefault}{\mddefault}{\updefault}{$\beta$}%
}}}
\put(3751,614){\makebox(0,0)[lb]{\smash{\SetFigFont{10}{12.0}{\familydefault}{\mddefault}{\updefault}{$\gamma$}%
}}}
\put(4501,-1936){\makebox(0,0)[lb]{\smash{\SetFigFont{10}{12.0}{\familydefault}{\mddefault}{\updefault}{$\delta$}%
}}}
\put(7951,-1111){\makebox(0,0)[lb]{\smash{\SetFigFont{10}{12.0}{\familydefault}{\mddefault}{\updefault}{$\delta$}%
}}}
\put(7876,1514){\makebox(0,0)[lb]{\smash{\SetFigFont{10}{12.0}{\familydefault}{\mddefault}{\updefault}{$\beta$}%
}}}
\put(7876,1064){\makebox(0,0)[lb]{\smash{\SetFigFont{10}{12.0}{\familydefault}{\mddefault}{\updefault}{$\gamma$}%
}}}
\put(7876,614){\makebox(0,0)[lb]{\smash{\SetFigFont{10}{12.0}{\familydefault}{\mddefault}{\updefault}{$\delta$}%
}}}
\put(7426,-1861){\makebox(0,0)[lb]{\smash{\SetFigFont{10}{12.0}{\familydefault}{\mddefault}{\updefault}{$\beta$}%
}}}
\put(7426,-2311){\makebox(0,0)[lb]{\smash{\SetFigFont{10}{12.0}{\familydefault}{\mddefault}{\updefault}{$\gamma$}%
}}}
\put(7426,-2761){\makebox(0,0)[lb]{\smash{\SetFigFont{10}{12.0}{\familydefault}{\mddefault}{\updefault}{$\delta$}%
}}}
\put(7051,-3736){\makebox(0,0)[lb]{\smash{\SetFigFont{10}{12.0}{\familydefault}{\mddefault}{\updefault}{$\beta$}%
}}}
\put(7951,-3736){\makebox(0,0)[lb]{\smash{\SetFigFont{10}{12.0}{\familydefault}{\mddefault}{\updefault}{$\gamma$}%
}}}
\put(7501,-4711){\makebox(0,0)[lb]{\smash{\SetFigFont{10}{12.0}{\familydefault}{\mddefault}{\updefault}{$\delta$}%
}}}
\put(601,-2986){\makebox(0,0)[lb]{\smash{\SetFigFont{10}{12.0}{\familydefault}{\mddefault}{\updefault}{$X_3$}%
}}}
\put(601,-1186){\makebox(0,0)[lb]{\smash{\SetFigFont{10}{12.0}{\familydefault}{\mddefault}{\updefault}{$X_2$}%
}}}
\put(601,614){\makebox(0,0)[lb]{\smash{\SetFigFont{10}{12.0}{\familydefault}{\mddefault}{\updefault}{$X_1$}%
}}}
\put(5401,1214){\makebox(0,0)[lb]{\smash{\SetFigFont{10}{12.0}{\familydefault}{\mddefault}{\updefault}{$Y_1$}%
}}}
\put(5401,-361){\makebox(0,0)[lb]{\smash{\SetFigFont{10}{12.0}{\familydefault}{\mddefault}{\updefault}{$Y_2$}%
}}}
\put(5401,-2161){\makebox(0,0)[lb]{\smash{\SetFigFont{10}{12.0}{\familydefault}{\mddefault}{\updefault}{$Y_3$}%
}}}
\put(5401,-3961){\makebox(0,0)[lb]{\smash{\SetFigFont{10}{12.0}{\familydefault}{\mddefault}{\updefault}{$Y_4$}%
}}}
\put(1426,-586){\makebox(0,0)[lb]{\smash{\SetFigFont{10}{12.0}{\familydefault}{\mddefault}{\updefault}{$\alpha$}%
}}}
\put(1351,-2386){\makebox(0,0)[lb]{\smash{\SetFigFont{10}{12.0}{\familydefault}{\mddefault}{\updefault}{$\alpha$}%
}}}
\put(1351,-4111){\makebox(0,0)[lb]{\smash{\SetFigFont{10}{12.0}{\familydefault}{\mddefault}{\updefault}{$\alpha$}%
}}}
\put(3976,-1936){\makebox(0,0)[lb]{\smash{\SetFigFont{10}{12.0}{\familydefault}{\mddefault}{\updefault}{$\gamma$}%
}}}
\put(6376,464){\makebox(0,0)[lb]{\smash{\SetFigFont{10}{12.0}{\familydefault}{\mddefault}{\updefault}{$\alpha$}%
}}}
\put(6376,-1186){\makebox(0,0)[lb]{\smash{\SetFigFont{10}{12.0}{\familydefault}{\mddefault}{\updefault}{$\alpha$}%
}}}
\put(6376,-3361){\makebox(0,0)[lb]{\smash{\SetFigFont{10}{12.0}{\familydefault}{\mddefault}{\updefault}{$\alpha$}%
}}}
\put(6301,-5161){\makebox(0,0)[lb]{\smash{\SetFigFont{10}{12.0}{\familydefault}{\mddefault}{\updefault}{$\alpha$}%
}}}
\put(7501,-136){\makebox(0,0)[lb]{\smash{\SetFigFont{10}{12.0}{\familydefault}{\mddefault}{\updefault}{$\beta$}%
}}}
\put(8401,-136){\makebox(0,0)[lb]{\smash{\SetFigFont{10}{12.0}{\familydefault}{\mddefault}{\updefault}{$\gamma$}%
}}}
\end{picture}
    \caption{Thickenings of $X$ and $Y$, always given
    by the immersion in the plane} \label{thick-graphs5-3:fig}
    \end{center}
    \end{figure}

Up to symmetry, the possible thickenings of $X$ and $Y$ to
orientable surfaces with boundary are those described in
\fullref{thick-graphs5-3:fig}, as explained in the caption.
The associated closed surfaces $\Sigmatil$ and the half-lengths of the
discs added are as follows:
\begin{eqnarray*}
X_1:&\ &\matS,\ (\alpha+1,\gamma+1,\delta+1,\alpha+2\beta+\gamma+\delta+5) \\
X_2:&\ &\matS,\ (\alpha+1,\gamma+1,\gamma+\delta+2,\alpha+2\beta+\delta+4) \\
X_3:&\ &\matT,\ (\alpha+1,\alpha+2\beta+2\gamma+2\delta+7) \\
Y_1:&\ &\matS,\ (\alpha+1,\beta+\gamma+2,\gamma+\delta+2,\alpha+\beta+\delta+3) \\
Y_2:&\ &\matT,\ (\alpha+1,\alpha+2\beta+2\gamma+2\delta+7) \\
Y_3:&\ &\matT,\ (\beta+\gamma+2,2\alpha+\beta+\gamma+2\delta+6) \\
Y_4:&\ &\matT,\
(\alpha+\beta+\gamma+3,\alpha+\beta+\gamma+2\delta+5) .
\end{eqnarray*}
We must then determine all the possible values which can be
attained by these strings as $\alpha,\beta,\gamma,\delta$ vary
among non-negative integers under the constraint
$\alpha+\beta+\gamma+\delta=d/2-4$.

Let us now specialize the proof for $\Sigmatil=\matT$,
corresponding to cases $X_3$, $Y_2$, $Y_3$, and $Y_4$. It is
obvious that the string $(d/2,d/2)$ cannot be realized. If
$1\leqslant k\leqslant d/2-3$ we can realize $(d/2+k,d/2-k)$ using
$Y_4$ with $\alpha=d/2-3-k$, $\beta=\gamma=0$, $\delta=k-1$. Using
$X_3$ and $Y_3$ we can obviously realize $(d-1,1)$ and $(d-2,2)$.
This proves the first assertion.

Turning to $\Sigmatil=\matS$, let us
first concentrate on case $Y_1$. We denote the unordered elements
of $(d_{3j})_{j=1,\ldots,4}$ by $x,y,z,w=d-x-y-z$, and solve in
the unknowns $\alpha,\beta,\gamma,\delta$ the system
$$\left\{
\begin{array}{rcl}
x & = & \alpha+1 \\
y & = & \beta+\gamma+2 \\
z & = & \gamma+\delta+2 \\
w & = & \alpha+\beta+\delta+3
\end{array}\right.
\quad\Rightarrow\quad \left\{
\begin{array}{rcl}
\alpha & = & x-1 \\
\beta & = & d/2-x-z-1 \\
\gamma & = & x+y+z-1-d/2 \\
\delta & = & d/2-x-y-1.
\end{array}\right.$$
This solution is acceptable if and only if
$\alpha,\beta,\gamma,\delta$ are all non-negative, namely if the
following holds:
$$(*)\qquad\left\{
\begin{array}{rcl}
x + y & \leqslant & d/2-1 \\
x + z & \leqslant & d/2-1 \\
x+y+z & \geqslant & d/2+1.
\end{array}\right.$$
We deduce that a branch datum with $\Sigmatil=\matS$ is realizable
using case $Y_1$ if and only if it is possible to extract from
$(d_{3j})_{j=1,\ldots,4}$ integers $x,y,z$ satisfying $(*)$. Let
us now denote by $\ell$ the largest of the $d_{3j}$'s and prove
the following facts:

\textbf{Claim 1}\qua \textsl{If $\ell\geqslant d/2$ then
the branch datum cannot be realized using $Y_1$.}
\vspace{10pt}

\textbf{Claim 2}\qua \textsl{If $\ell< d/2$ then the branch
datum can be realized using $Y_1$ if and only if it does not have
the form $(k,k,d/2-k,d/2-k)$.}
\vspace{10pt}

Claim 1 is easy: of course we cannot choose $x$, $y$, or $z$ to be
$\ell$, otherwise one of the first two conditions in $(*)$ would
be violated, whence $x+y+z=d-\ell\leqslant d-d/2=d/2$, which
contradicts the last condition in $(*)$. Turning to Claim 2, it is
clear that from $(k,k,d/2-k,d/2-k)$ we cannot extract $x,y,z$
satisfying $(*)$. To prove the converse, let us choose $x,y,z,w$
in increasing order, ie,
$$x\leqslant y\leqslant z\leqslant w=\ell=d-x-y-z<d/2$$
which implies the last condition in $(*)$. Under these assumptions
the second condition in $(*)$ implies the first one. If the second
condition is violated then we have
\begin{eqnarray*}
& & x\leqslant y\leqslant z=:d/2-t\leqslant \ell=:d/2-s \\
& & x+z\geqslant d/2 \Rightarrow x\geqslant t \\
& & x+y=t+s\Rightarrow y\leqslant s \\
& & z\leqslant\ell \Rightarrow t\geqslant s.
\end{eqnarray*}
These facts imply that $x\leqslant y\leqslant s\leqslant
t\leqslant x$. Calling $k$ this common value, we deduce that
$(x,y,z,w)$ has the form $(k,k,d/2-k,d/2-k)$, and the claim is
established.
\vspace{10pt}

To conclude the proof for
$\Sigmatil=\matS$ it is now sufficient to establish the following:
\vspace{10pt}

\textbf{Claim 3}\qua \textsl{If $\ell\geqslant d/2$ then a
branch datum can be realized using $X_1$ or $X_2$ if and only if
$(d_{3j})_{j=1,\ldots,4}$ is not $(d/2,d/6,d/6,d/6)$.}
\vspace{10pt}

\textbf{Claim 4}\qua \textsl{No branch datum with
$(d_{3j})_{j=1,\ldots,4}$ of the form $(k,k,d/2-k,d/2-k)$ can be
realized using $X_1$ or $X_2$.}
\vspace{10pt}

Let us prove Claim 3. Again we denote the $d_{3j}$'s by $x,y,z,w$,
and we choose them in increasing order, ie,
$$x\leqslant y\leqslant z<w=\ell=d-x-y-z,$$
whence $x+y+z\leqslant d/2$. If we want to realize the datum using
$X_1$ we have the following forced choice up to symmetry:

$$\left\{
\begin{array}{rcl}
x & = & \alpha+1 \\
y & = & \gamma+1 \\
z & = & \delta+1 \\
w & = & \alpha+2\beta+\gamma+\delta+5
\end{array}\right.
\qquad\Rightarrow\qquad \left\{
\begin{array}{rcl}
\alpha & = & x-1\\
\beta & = & d/2-x-y-z-1 \\
\gamma & = & y-1\\
\delta & = & z-1.
\end{array}\right.$$

This is an acceptable solution unless $x+y+z=d/2$. So we turn to
case $X_2$ and try to realize the case $x\leqslant y\leqslant
z<\ell=x+y+z$. We first solve

$$\left\{
\begin{array}{rcl}
x & = & \alpha+1 \\
y & = & \gamma+1 \\
z & = & \gamma+\delta+2 \\
x+y+z & = & \alpha+2\beta+\delta+4
\end{array}\right.
\qquad\Rightarrow\qquad \left\{
\begin{array}{rcl}
\alpha & = & x-1\\
\beta & = & y-1 \\
\gamma & = & y-1\\
\delta & = & z-y-1
\end{array}\right.$$

which is acceptable only for $z>y$. We are left to deal with the
case $x\leqslant y=z<\ell=x+2y$ and we try to use $X_2$ solving

$$\left\{
\begin{array}{rcl}
y & = & \alpha+1 \\
x & = & \gamma+1 \\
y & = & \gamma+\delta+2 \\
x+2y & = & \alpha+2\beta+\delta+4
\end{array}\right.
\qquad\Rightarrow\qquad \left\{
\begin{array}{rcl}
\alpha & = & y-1\\
\beta & = & x-1 \\
\gamma & = & x-1\\
\delta & = & y-x-1
\end{array}\right.$$
which is acceptable only for $y>x$. We have thus realized the
branch datum for all $(d_{3j})_{j=1,\ldots,4}$ except
$(d/2,d/6,d/6,d/6)$, which of course cannot be realized.

Turning to Claim 4, we begin using $X_2$. Since
$\gamma+1<\gamma+\delta+2$ and $\alpha+1<\alpha+2\beta+\delta+4$,
there is only one attempt we can make:
$$\left\{
\begin{array}{rcl}
k & = & \alpha+1 \\
k & = & \gamma+1 \\
d/2-k & = & \gamma+\delta+2 \\
d/2-k & = & \alpha+2\beta+\delta+4
\end{array}\right.
\qquad\Rightarrow\qquad \left\{
\begin{array}{rcl}
\alpha & = & k-1\\
\beta & = & -1 \\
\gamma & = & k-1\\
\delta & = & d/2-2k-1
\end{array}\right.$$
which is not acceptable. We then try $X_1$. Up to symmetry there
is again one case only, which is of course impossible. \finedimo
\vspace{5pt}

A similar argument allows one to prove
\fullref{23:nonex:prop}.
\vspace{5pt}

\paragraph{Exceptional data with non-prime degree} Here we obtain
the proof of \fullref{excep:by:fixpoints:thm}. We start with a
lemma.
\vspace{5pt}

\begin{lemma}\label{transpos-excep:lem}
Suppose that $d=kh$ with $k,h\geqslant 2$. Let
$(s_j)_{j=1,\ldots,p}$, $(t_j)_{j=1,\ldots,q}$ be partitions
of $h$ with $p,q\geqslant 2$ and $p+q\geqslant h+2$. Then for all $1\leqslant r<(p+q-h)$,
if $n=p+q-r-h+2$, the following branch datum is exceptional:
$$\begin{array}{r}
\big(\matS,\matS,n,d,(ks_1,\ldots,ks_p),(kt_1,\ldots,kt_q),(h+r,1,\ldots,1), \\
(2,1,\ldots,1),\ldots,(2,1,\ldots,1)\big).
\end{array}$$
\end{lemma}
\vspace{5pt}

\begin{proof}
It is easy to see that the datum is compatible. The proof of
exceptionality is by induction on $n$. The base of induction is
$n=3$. We remark that in this case the statement of the lemma
could be inferred from results announced by Edmonds, Kulkarni, and
Stong~\cite[page~775]{EKS}, which are however stated without
proof, so for the sake of completeness we provide an independent
argument.
\vspace{5pt}

When $n=3$, the branch datum does not contain partitions
$(2,1,\ldots,1)$. To prove the lemma in this case, we proceed by
induction on $h$. If $h=2$ we have $p=q=2$ and
$(t_1,t_2)=(s_1,s_2)=(1,1)$, while for $h\geqslant 3$ we have two
cases: either one of the partitions $(t_j)_{j=1,\ldots,p}$ and
$(s_j)_{j=1,\ldots,q}$ has the form $(1,\ldots,1)$ or not. So the
next claim serves both to prove the base step of the induction and
to deal with the first case of the inductive step.
\vspace{5pt}

\textbf{Claim}\qua \textsl{Suppose that $k,h\geqslant 2$. Let
$(s_j)_{j=1,\ldots,p}$ be a partition of $h$ with $p\geqslant 2$.
Then the branch datum
$\big(\matS,\matS,3,kh,(ks_1,\ldots,ks_p),(k,\ldots,k),
(h+p-1,1,\ldots,1)\big)$ is exceptional.} 
\vspace{5pt}

First note that the number of $1$'s in the third partition is
$(k-1)h-p+1$. Take arbitrary $\tau_1,\tau_2\in\permu_{kh}$ with cyclic
structures $(ks_1,\ldots,ks_p)$ and $(k,\ldots,k)$ such that $\langle
\tau_1,\tau_2\rangle$ is transitive, and consider the associated graph
$D(\tau_1,\tau_2)$. If $\tau_1$, $\tau_2$ realize the above datum,
$\tau_1\cdot\tau_2$ has $(k-1)h-p+1$ fixed points. Remark that to each
fixed point of $\tau_1\cdot\tau_2$ there corresponds in
$D(\tau_1,\tau_2)$ a pair of edges having the same ends. More
precisely, let us define a multi-edge of $D(\tau_1,\tau_2)$ as the set
of all edges having two given vertices as ends. Then one sees that a
multi-edge $\varphi$ gives rise to at most $\#(\varphi)$ fixed points
of $\tau_1\cdot\tau_2$, and all fixed points arise like this. The case
of $\#(\varphi)$ fixed points actually occurs only if $\tau_1$ and
$\tau_2$ contain cycles which are the inverse of each other, but this
is easily recognized to be incompatible with the assumptions of the
claim. So a multi-edge $\varphi$ contributes with at most
$\#(\varphi)-1$ fixed points. We conclude that the total number of
fixed points is at most the sum of all the multiplicities of the
multi-edges, which is equal to $kh$, minus the number of multi-edges.
\vspace{5pt}

Let us now estimate the number of multi-edges.  By definition of
$D(\tau_1,\tau_2)$, the set of its vertices is split as $V_1\sqcup
V_2$, where $V_1$ consists of $p$ vertices of valences
$ks_1,\ldots,ks_p$ and  $V_2$ consists of $h$ vertices of valence
$k$. This implies that at least $s_j$ multi-edges are incident to
the $j$-th vertex of $V_1$. However, if there are exactly $s_j$
multi-edges, connectedness of $D(\tau_1,\tau_2)$ easily implies
that $p=1$, which was excluded. So there are at least $s_j+1$
multi-edges. Therefore there are in total at least
$\sum_{j=1}^p(s_j+1)=h+p$ multi-edges. From the above we deduce
that $\tau_1\cdot\tau_2$ has at most $kh-(h+p)=(k-1)h-p$ fixed
points. This implies that $\tau_1,\tau_2$ cannot realize the
branch datum, hence the Claim is established.
\vspace{5pt}

To conclude the inductive proof on $h$ (still with $n=3$, which is
the base of our more general induction) we have to deal with the
case where $h\geqslant 3$ and both partitions
$(s_j)_{j=1,\ldots,p}$ and $(t_j)_{j=1,\ldots,q}$ contain at least
one entry larger than 1. This implies that $p,q\leqslant h-1$ and,
taking into account the inequality $p+q\geqslant h+2$, that
$p,q\geqslant 3$. Suppose that a datum as in the statement of the
lemma (under the current assumptions) is realizable and consider
the dessin d'enfant $D$ constructed as in the proof of
\fullref{dessins:coverings:prop}.
\newpage

We first show that \emph{$D$ contains vertices $v$ and $u$ such
that all the edges of $D$ incident to $v$ join $v$ to $u$.}
Indeed, notice that $\matS^2\setminus D$ consists of some bigons
and one $2(p+q-1)$--gon. Successively compressing each bigon into a
single edge (see \fullref{multi:edge:fig}),
\begin{figure}[ht!]
    \begin{center}
    \includegraphics[scale=0.5]{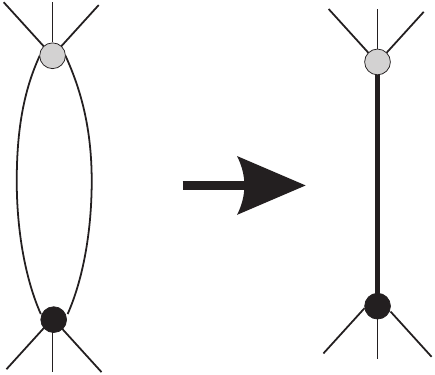}
    \caption{The passage from $D$ to $\overline{D}$}
    \label{multi:edge:fig}
    \end{center}
    \end{figure}
we get an embedded graph $\overline{D}$ in $\matS^2$ with the same
vertices as $D$, whose complement consists of a single disc.
Therefore $\overline{D}$ is a tree and not a point, so it contains
vertices of valence 1. Any such vertex $v$ and the vertex $u$
joined to it in $\overline{D}$ has the desired property.
\vspace{5pt}

Now let $D'$ be the graph obtained from $D$ by deleting $v$ and
all the edges incident to it. Obviously, $D'$ is still a dessin
d'enfant and its complement consists of some bigons and one
$2(p+q-2)$--gon. Without loss of generality, we may assume that
$v\in V_2$ and $u\in V_1$. If the valence of $v$ in $D$ is $kt_a$
and that of $u$ is $ks_b$ then the valence of $u$ in $D'$ is
$k(s_b-t_a)$ and this number is positive, otherwise $p=q=1$, which
is excluded. Therefore $D'$ realizes a branch datum as in the
statement, with $n=3$, $k'=k$, $h'=h-t_a$, $p'=p$, $q'=q-1$, a
partition $(s'_j)$ of $h'$ obtained from $(s_j)$ by replacing
$s_b$ by $s_b-t_a$ and reordering, and another partition $(t'_j)$
of $h'$ obtained from $(t_j)$ by dropping $t_a$. To conclude the
proof for the case $n=3$ we must show that the conditions
$k',h',p',q'\geqslant 2$ and $p'+q'\geqslant h'+2$ are fulfilled.
Of course $k',p'\geqslant 2$. Moreover $q'\geqslant 2$ because
$q\geqslant 3$, which implies that $h'\geqslant 2$. Now
$h'=h-t_a\leqslant h-1$ and $p+q\geqslant h+2$, so
$$p'+q'=p+q-1\geqslant h+2-1=h+1\geqslant h'+2.$$
The case $n=3$ is eventually settled.
\vspace{5pt}

Suppose now that $n\geqslant 4$, so the last partition has the
form $(2,1,\ldots,1)$. This implies that $p+q-h\geqslant
r+2\geqslant 3$, and, since $p,q\leqslant h$, we have
$p,q\geqslant 3$. Assume, by contradiction, that the datum is
realizable by a map $f$. Let $y_i$ be the branching point
corresponding to the $i$-th partition and consider the dessin
d'enfant $D$ constructed as in the proof of
\fullref{dessins:coverings:prop} with $x_{n-2}=y_2$,
$x_{n-1}=y_n$, $x_n=y_1$ (the other $x_j$'s must be the other
$y_j$'s, the order does not matter). In particular $D$ contains
$f^{-1}(\alpha_{n-2}\cup\alpha_{n-3})$, where $\alpha_{n-2}$ and
$\alpha_{n-3}$ join $y_n$ to $y_2$ and $y_2$ to some other $y_j$
respectively. Recall that $V_i$ is $f^{-1}(x_i)$, so
$V_{n-1}=f^{-1}(y_n)$ and $V_{n-2}=f^{-1}(y_2)$. Considering
valences, one concludes that there are only two possibilities for
$f^{-1}(\alpha_{n-2}\cup\alpha'_{n-3})$, where $\alpha'_{n-3}$ is
the half of $\alpha_{n-3}$ incident to $y_2$. These possibilities
are shown in \fullref{transpos:fig}, where the elements of
$V_{n-1}$ are the black dots and those of $V_{n-2}$ are the grey
dots.
\begin{figure}[ht!]
    \begin{center}
    \includegraphics[scale=0.5]{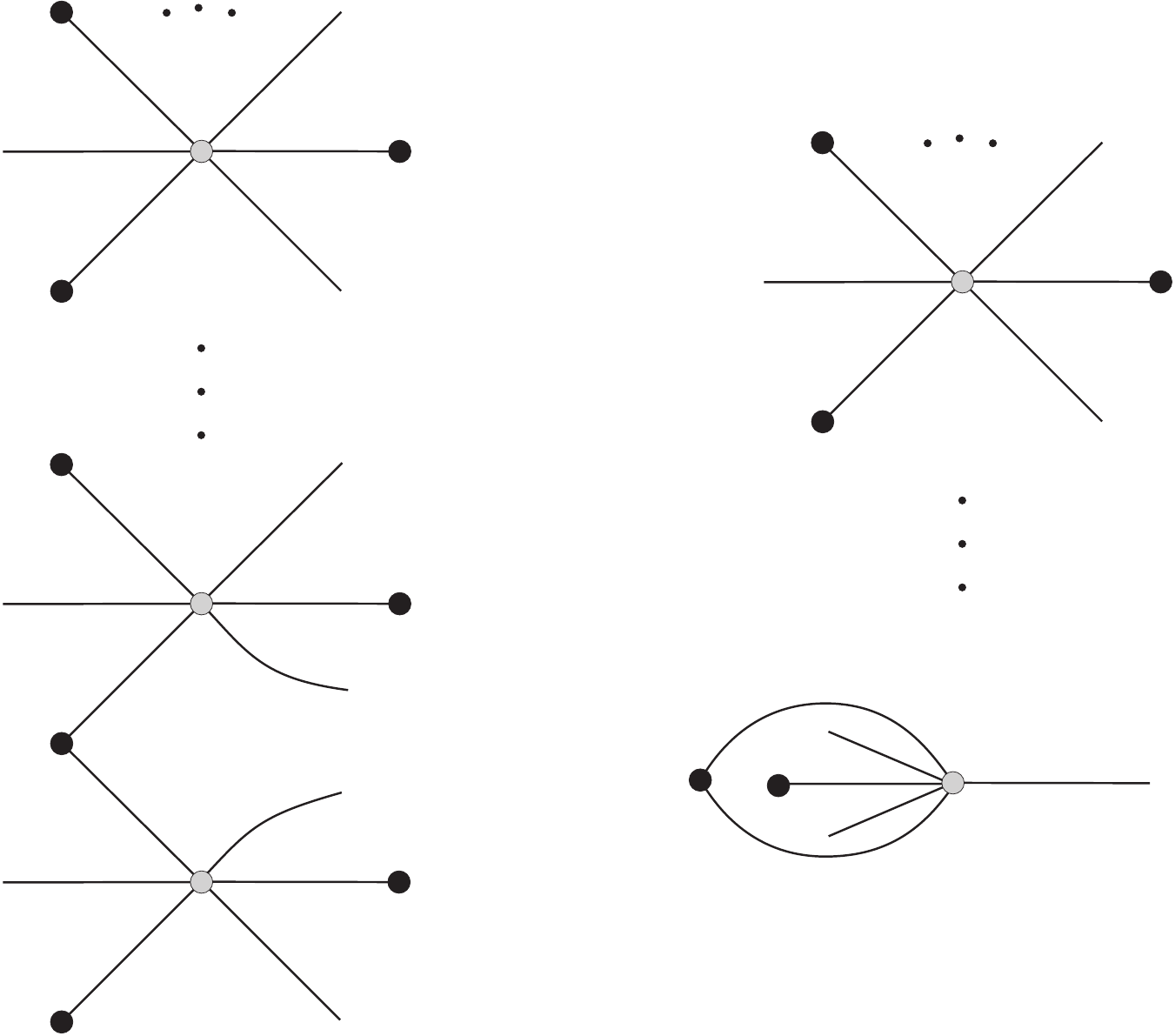}
    \caption{The two possibilities for $f^{-1}(\alpha_{n-2}\cup\alpha'_{n-3})$}
    \label{transpos:fig}
    \end{center}
    \end{figure}
Let $e_{n-2}^{(i)}$, $e_{n-2}^{(j)}$ be the edges with an endpoint
at the 2--valent black vertex. Suppose we have the case shown in
\fullref{transpos:fig}-left, ie, that $e_{n-2}^{(i)}$ and
$e_{n-2}^{(j)}$ are incident to two distinct grey vertices, which
we denote by $v'$ and $v''$. Then we remove all univalent black
vertices together with the edge adjacent to any such vertex, and
contract the set $e_{n-2}^{(i)}\cup e_{n-2}^{(j)}$ to a point,
which now becomes a new grey vertex. This gives a new dessin
d'enfant $D'$ and we can now analyze which branch datum it
realizes. We first notice that there is a natural correspondence
between the complementary discs of $D'$ and those of $D$, and the
length of the $j$-th one of $D'$ is $2d_{1j}=2ks_j$ less than that
of the $j$-th one of $D$. Moreover the contraction leading from
$D$ to $D'$ fuses $v'$ and $v''$ together, so the number of grey
vertices is decreased by one. Since exactly half of all the edges
incident to a grey vertex belong to $E_{n-2}$, and all these edges
are destroyed (either by removal or by contraction), the valence
of any grey vertex distinct from $v'$, $v''$ is halved, and the
valence of the new vertex obtained by fusing is half the sum of
the valences of $v'$ and $v''$. However the grey vertices were
non-terminal in $D$ but they become terminal in $D'$. All this
shows that $D'$ is a dessin d'enfant realizing, up to re-ordering
the $t_j$'s and the entries in the second partition, the datum
$\big(\matS,\matS,n-1,d,(ks_1,\ldots,ks_p),(k(t_1+t_2),kt_3,\ldots,kt_q),(h+r,1,\ldots,1),
(2,1,\ldots,1),\ldots,(2,1,\ldots,1)\big)$. Since $q\geqslant 3$
we again get a datum of the form described in the statement, but
with $n-1$ in place of $n$. This contradicts the inductive
assumption.
\newpage

If we have the case shown in \fullref{transpos:fig}-right, we
remove all the black vertices and all the edges incident to them.
Then the valence of each grey vertex gets halved, and two of the
complementary discs are merged into a single disc. A discussion
similar to the above one shows that we get a dessin d'enfant
realizing (up to re-ordering) the datum
$\big(\matS,\matS,n-1,d,(k(s_1+s_2),ks_3,\ldots,ks_p),(kt_1,\ldots,kt_q),(h+r,1,\ldots,1),
(2,1,\ldots,1),\ldots,(2,1,\ldots,1)\big)$. As before, the datum
is non-realizable by the inductive assumption, whence the
conclusion in all cases.
\end{proof}

\dimo{excep:by:fixpoints:thm} We assume $n\geqslant 3$ otherwise
the statement is empty. We first claim that $m_1,m_2\geqslant 2$.
By contradiction, suppose that $d_{11}=d$, all the $d_{2j}$'s are
multiples of $k$, and $d_{31}>d/k$. Then $m_1=1$, $m_2\leqslant
d/k$, $m_3\leqslant d-d/k$, and $m_i\leqslant d-1$ for $i\geqslant
4$. By the Riemann--Hurwitz condition
$(n-2)d+2=m_1+\ldots+m_n\leqslant 1+d/k+(d-d/k)+(n-3)(d-1)$, which
implies that $n\leqslant 2$, a contradiction.

Suppose now that there is a realizable datum of the form
$$\big(\matS,\matS,n,d,(ks_1,\ldots,ks_p),(kt_1,
\ldots,kt_q),(d/k+r,d_{32},\ldots,d_{3m_3}),
(d_{ij})_{j=1,\ldots,m_i}\big),$$ with $r\geqslant 1$ and $p,q\geqslant 2$. By
\fullref{Hurwitz:method:thm} we can find permutations
$\tau_1$, $\ldots$, $\tau_n$ such that the cyclic structure of
$\tau_i$ is given by the $i$-th partition and
$\tau_1\tau_2\ldots\tau_n=1$. Obviously, we can present $\tau_3$
as a product of a cycle $\tau_3'$ of length $d/k+r$ and a
permutation $\tau_3''$ of cyclic structure
$(1,\ldots,1,d_{32},\ldots,d_{3m_3})$. Notice that
$\tau_3''\tau_4\ldots\tau_n$ can be obtained as the product of
$$(d_{32}+\ldots+d_{3m_3}-m_3+1)+\sum_{i=4}^n(d-m_i)=p+q-r-d/k-1$$
transpositions. The collection of these transpositions together
with the permutations $\tau_1,\tau_2,\tau_3'$ provides thus a
realization of the datum
$$\begin{array}{r}
\big(\matS,\matS,p+q-r-d/k+2,d,(ks_1,\ldots,ks_p),(kt_1,\ldots,kt_q),(d/k+r,1,\ldots,1), \\
(2,1,\ldots,1),\ldots,(2,1,\ldots,1)\big),
\end{array}$$
which is non-realizable by \fullref{transpos-excep:lem}, whence
the conclusion. \finedimo

\section{Exceptions due to decomposability}\label{nonex:decom:section}

We now describe another technique that can be employed to prove
exceptionality of branch data. The basic underlying remark is that
certain patterns in a branch datum force the covering realizing
the datum to be \emph{decomposable}, namely the composition of two
non-trivial coverings, and one can deduce strong restrictions on
the whole datum from the information that all its realizations are
decomposable.

\paragraph{Decomposability of coverings and permutations}
We will now state an easy result characterizing decomposable
coverings, apparently due to Ritt and cited (in a less detailed
fashion) in~\cite[Theorem 1.7.6]{LandoZv}. We first introduce some
terminology and establish a lemma which are necessary for the
statement.

\begin{defn}
If $1<k<d$ and $k|d$ we call {\em block decomposition} of
{\em order} $k$ for $\tau\in\permu_d$ a partition of
$\{1,2,...,d\}$ into $d/k$ subsets (the {\em blocks}) such that
each block has $k$ elements and $\tau$ induces a well-defined
permutation $\widehat\tau$ of the blocks. We say that $G<\permu_d$
has a block decomposition of order $k$ if there is a common such
decomposition for all the elements of $G$.
\end{defn}

\begin{lemma}\label{exist:dec:lem}
Let $\tau$ have cyclic structure $(d_1,\ldots,d_m)$. Then $\tau$
admits a block decomposition of order $k$ if and only if
$(d_1,\ldots,d_m)$ can be partitioned into sets $D_1,\ldots,D_t$
and there exist integers $p_1,\ldots,p_t\geqslant 1$ such that for
all $j$
$$p_j|x\quad \forall x\in D_j
\qquad{\it and}\qquad \sum_{x\in D_j} \frac x{p_j}=k.$$ In this
case $\widehat\tau$ has cyclic structure $(p_1,\ldots,p_t)$.
\end{lemma}

\begin{proof}
Suppose the block decomposition exists, and take a cycle
$(B_1,\ldots,B_p)$ of $\widehat\tau$, with
$B_{j}=\{r_{j,1},\dots,r_{j,k}\}$. Since $\tau(B_j)=B_{j+1}$ for
$j=1,\ldots,p-1$, up to changing notation we have
$$\tau(r_{1,1})=r_{2,1}\qquad \tau(r_{2,1})=r_{3,1}\qquad\ldots\qquad
\tau(r_{p-1,1})=r_{p,1}.$$ Now we either have
$\tau(r_{p,1})=r_{1,1}$ or, up to changing notation,
$\tau(r_{p,1})=r_{1,2}$. Proceeding similarly we see that $\tau$
contains a cycle
$$(r_{1,1},\dots,r_{p,1},r_{1,2},\ldots,r_{p,2},\ldots,r_{1,b_1},\ldots,r_{p,b_1})$$
of order $x_1=p\cdot b_1$. Repeating the same argument from
$r_{1,b_1+1}$ and so on we find cycles of length $x_1,\ldots,x_s$
of $\tau$ such that $p|x_l$ and $\sum_{l=1}^{s}\frac {x_l}p=k$.
The conclusion follows by considering the totality of the cycles
of $\widehat\tau$, and the opposite implication is proved along
the same lines.
\end{proof}

\begin{cor}\label{Ritt:cor}
Let $\theta\co\pi_1(\Sigma_n)\to\permu_d$ be associated as in
\fullref{Hurwitz:method:thm} to a branched covering realizing
a branch datum $\big(\Sigmatil,\Sigma,n,d,(d_{ij})\big)$. Then the
covering is decomposable if and only if $\mbox{Im}(\theta)$ has a
block decomposition. Moreover, if $\mbox{Im}(\theta)$ has a block
decomposition of order $k$ with the decomposition of $\theta(e_i)$
corresponding to a partition of $(d_{ij})$ into
$D_{i1},\ldots,D_{it_i}$ and integers $p_{i1},\ldots,p_{it_i}$
then the covering factors through coverings realizing branch data
\begin{eqnarray*}
& \left(\Sigmatil,\Sigma',t_1+\ldots+t_n,k,\left(\left(\frac
x{p_{ij}}\right)_{x\in
D_{ij}}\right)_{i=1,\ldots,n,\ j=1,\ldots,t_i}\right) & \\
&
\left(\Sigma',\Sigma,n,d/k,\left(\left(p_{ij}\right)_{j=1,\ldots,t_i}\right)_{i=1,\ldots,n}\right).
&
\end{eqnarray*}
\end{cor}

Note that to apply this result for any given $\theta$ one could
use~\cite{Atkinson}.

\paragraph{Very even data}
Here we apply the methods of the previous paragraph and those of
\fullref{nonex:dessin:section} to describe two more infinite
classes of exceptional data. \fullref{even-deg_exceptions:thm}
is in fact an immediate consequence of the following result:

\begin{prop}\label{filtration:prop}
If $d$ is even and all $d_{ij}$ are also even for $i=1,2$ then any
covering realizing a branch datum
$\big(\matS,\matS,n,d,(d_{ij})\big)$,
up reordering the $d_{ij}$'s, is a composition of coverings
realizing data of the form
$$\begin{array}{rl}
\big(\matS,\matS,2n-2,d/2, & (d_{1j}/2)_{j=1,\ldots,m_1},(d_{2j}/2)_{j=1,\ldots,m_2},\\
& (d_{3j})_{j=1,\ldots,h_3},(d_{3j})_{j=h_3+1,\ldots,m_3},\ldots,\\
& (d_{nj})_{j=1,\ldots,h_n},(d_{nj})_{j=h_n+1,\ldots,m_n}\big),
\end{array}
$$
with $1\leqslant h_i<m_i$ for $i=3,\ldots,n$, and
$$\big(\matS,\matS,n,2,(2),(2),(1,1),\ldots,(1,1)\big).$$
\end{prop}

The proof requires a definition and an easy lemma. We call
\emph{checkerboard graph} a finite $1$--subcomplex of a surface
whose complement consists of open discs each bearing a color black
or white, so that each edge separates black from white.

\begin{lemma}\label{even:check:lem}
A connected graph in $\matS$ with all vertices of even valence is
a checkerboard graph.
\end{lemma}

\begin{proof}
By induction on the number $v$ of vertices of the graph $D$. If
$v=1$ then $D$ is a wedge of circles, and the conclusion follows
from the fact that embedded circles are separating on $\matS$. If
$v>1$ choose an edge of $D$ having distinct ends, and let $D'$ be
obtained from $D$ by contracting this edge to a vertex. Then the
hypothesis applies to $D'$, so there is a checkerboard coloring of
$\matS\setminus D'$. Now the regions of $\matS\setminus D$ can be
naturally identified to those of $\matS\setminus D'$, and it is
easy to see that the coloring of $\matS\setminus D'$ also works
for $\matS\setminus D$.
\end{proof}

\dimo{filtration:prop}
Let $\tau_1,\ldots,\tau_{n-1}\in\permu_d$ correspond as in
\fullref{permu:S:n:cor} to a realization of the
given branch datum, with indices arranged so that
$\tau_1$ and $\tau_{n-1}$ correspond to the first two partitions of $d$.
We will prove that $\langle\tau_1,\ldots,\tau_{n-1}\rangle$ has a block
decomposition of order $d/2$ and then apply \fullref{Ritt:cor}.
Let $D=D(\tau_1,\ldots,\tau_{n-1})$ be the dessin d'enfant
associated to $\tau_1,\ldots,\tau_{n-1}$ as in
\fullref{permu:to:ext:dessin:prop}. Then $D$ is a
connected graph with all vertices of even valence, hence by
\fullref{even:check:lem} the components of $\matS\setminus D$
have a checkerboard coloring.

We will now choose black and white colors also for the edges of
$D$. We begin with the edges in $E_1$ and color them so that each
$V_1$--corner of a black component of $\matS\setminus D$ has first
a black and then a white edge in the positive order around the
vertex. Recalling that $\matS\setminus D$ is checkerboard colored
and that around each $v\in V_2$ the edges in $E_1$ and in $E_2$
alternate with each other, it is easy to see that all the edges in
$E_1$ incident to $v$ have the same color. We then color the edges
in $E_2$ incident to $v$ with the other color. Repeating this for
all $v$ we color $E_2$, and we can proceed similarly to color all
the edges of $D$.

Notice that by construction around any vertex of $D$ the colors of
edges alternate, therefore for $i\geqslant 2$ the edges
$e_{i-1}^{(k)}$ and $e_i^{(k)}$ always have opposite colors. Thus,
if we define
$$B_i=\{k:\ e_i^{(k)}\ {\rm is\ black}\},\qquad W_i=\{k:\
e_i^{(k)}\ {\rm is\ white}\}$$ then we have $B_1=W_2=B_3=\ldots$
and $W_1=B_2=W_3=\ldots$. We denote the former set by $B$ and the
latter set by $W$. The construction of $D$ given to prove
\fullref{permu:to:ext:dessin:prop} shows that
$\tau_2,\ldots,\tau_{n-2}$ leave $B$ and $W$ invariant, while
$\tau_1$ and $\tau_{n-1}$ switch them. Therefore we have a block
decomposition of $\langle\tau_1,\ldots,\tau_{n-1}\rangle$ of order
$d/2$. \fullref{exist:dec:lem} and \fullref{Ritt:cor} now
imply that the covering $\Sigmatil\to\matS$ realizing the given
branch datum factors through coverings $\Sigmatil\to\Sigma'$ and
$\Sigma'\to\Sigma$, with the latter having branch datum of the
form $\big(\matS,\matS,n,2,(2),(2),(1,1),\ldots,(1,1),*\big)$. But
the only such covering is as in the statement and the conclusion
easily follows.\finedimo

\dimo{mixed:deg:excep:cor} By \fullref{filtration:prop}
any covering realizing a datum as in the statement is a
composition of the appropriate covering of degree $2$ and a
covering given by the datum
$$\begin{array}{rl}
\big(\matS,\matS,2n-2,d/2, & (d_{11},\ldots,d_{1h_1}),
(d_{1(h_1+1)},\ldots,d_{1m_1}),\\
& (d_{21}/2,\ldots,d_{2m_2}/2),(d_{31}/2,\ldots,d_{3m_3}/2),\\
& (d_{4j})_{j=1,\ldots,h_4},(d_{4j})_{j=h_4+1,\ldots,m_4},\ldots,\\
& (d_{nj})_{j=1,\ldots,h_n},(d_{nj})_{j=h_n+1,\ldots,m_n}\big).
\end{array}
$$
So by \fullref{excep:by:fixpoints:thm} we have
$d_{ij}/2\leqslant d/2k$ for $i=2,3$ and $d_{ij}\leqslant d/2k$
for $i=4,\ldots,n$, whence the conclusion.\finedimo

We conclude the paper by establishing our only existence result.

\dimo{divisible:exist:thm} Let $\Sigma'$ be the orientable surface
with $\chi(\Sigma')=3-p$. Then the branch datum
$\big(\Sigma',\matS,3,p,(p),(p),(p)\big)$ is compatible, whence
realizable by \fullref{full:cycle:thm}. It is now easy to see
that the datum $\big(\Sigmatil,\Sigma',3,d/p,(d_{ij}/p)\big)$ is
also compatible. Now $\chi(\Sigma')\leqslant 0$, so the datum is
realizable by \fullref{OO:thm}, and the desired covering
$\Sigmatil\to\matS$ can then be constructed as the composition of
the coverings $\Sigmatil\to\Sigma'$ and
$\Sigma'\to\matS$.\finedimo

\bibliographystyle{gtart}
\bibliography{link}

\end{document}